\documentclass[smallextended,envcountsect,referee]{svjour3}

\smartqed
\usepackage{graphicx}
\journalname{JOTA}

\usepackage{amsmath}
\usepackage{mathrsfs}
\usepackage{setspace}
\usepackage{nicefrac}
\usepackage{graphicx}
\usepackage{float}
\usepackage{overpic}
\usepackage{amssymb}
\usepackage{color}
\usepackage{subcaption}
\graphicspath{{Figures/}}

\newenvironment{psmallmatrix}
  {\left(\begin{smallmatrix}}
  {\end{smallmatrix}\right)}

\newcommand{\sign}[0]{\text{sign}}
\usepackage{bm}

\makeatletter
\renewcommand{\maketag@@@}[1]{\hbox{\m@th\normalsize\normalfont#1}}%
\makeatother

\definecolor{yellow_color}{rgb}{0.9290, 0.6940, 0.1250}
\definecolor{green_color}{rgb}{0.4660, 0.6740, 0.1880}
\usepackage{bm}

\usepackage{cleveref}
\Crefname{equation}{}{Eqs.}
\Crefname{figure}{Figure}{Figures}
\Crefname{tabular}{Table}{Tables}

\definecolor{yellow_color}{rgb}{0.9290, 0.6940, 0.1250}
\definecolor{green_color}{rgb}{0.4660, 0.6740, 0.1880}

\begin{document}

\title{The Synthesis of Optimal Control Laws Using Isaacs' Method for the Solution of Differential Games}

\author{Meir Pachter and Isaac E. Weintraub}

\institute{Meir Pachter \at
             Air Force Institute of Technology \\
            Wright-Patterson AFB, OH 45433\\
            meir.pachter@afit.edu
           \and
            Isaac E. Weintraub  \at
            Air Force Research Laboratory \\
            Wright-Patterson AFB, OH 45433\\
            isaac.weintraub.1@us.af.mil
}

\date{Received: date / Accepted: date}

\maketitle

\begin{abstract}
In this paper we advocate for  Isaacs' method for the solution of differential games to be applied to the solution of optimal control problems. To make the argument, the vehicle employed is Pontryagin's canonical optimal control example, which entails a double integrator plant. However, rather than controlling the state to the origin, we require the end state to reach a terminal set that contains the origin in its interior. Indeed, in practice, it is required to control to a prescribed tolerance rather than reach a desired end state; constraining the end state to a terminal manifold of co-dimension n-1 renders the optimal control problem easier to solve. The global solution of the optimal control problem is obtained and the synthesized optimal control law is in state feedback form. In this respect, two target sets are considered: a smooth circular target and a square target with corners. Closed-loop state-feedback control laws are synthesized that drive the double integrator plant from an arbitrary initial state to the target set in minimum time. This is accomplished using Isaacs' method for the solution of differential games, which entails Dynamic Programming (DP), working backward from the Usable Part (UP) of the target set, as opposed to obtaining the optimal trajectories using the necessary conditions for optimality provided by Pontryagin's Maximum Principle (PMP). In this paper, the case is made for Isaacs' method for the solution of differential games to be applied to the solution of optimal control problems by way of the juxtaposition of the PMP and DP methods.
\end{abstract}



\section{Introduction}
In this paper, the Pontryagin Maximum Principle (PMP) and Dynamic Programming (DP) methods for the solution of optimal control problems are juxtaposed. We advocate for  Isaacs' method for the solution of differential games to be applied to the solution of optimal control problems.
The canonical example from \cite[pp. 23-27]{Pontryagin1962} concerning the application of the PMP to the synthesis of optimal controls is the vehicle employed to make the argument. The objective in \cite{Pontryagin1962} was to show that the necessary conditions for optimality embodied in the PMP yield a closed set of conditions such that the optimal control time-history can be obtained. The application of the PMP assumes the existence of an optimal control time-history and requires a Two-Point Boundary-Value Problem (TBVP) be solved, but with the provision that hard control constraints are allowed -- it is a necessary condition for optimality, akin to the situation in the calculus of variations. The objective of this work is to use the canonical example from reference \cite{Pontryagin1962} to demonstrate the application of differential game theory / Isaacs' method  \cite{Isaacs1965Differential} to optimal control problems and obtain their \emph{global} solution and state feedback optimal control laws. Isaacs' method is based on the constructive method of DP which provides sufficient conditions for optimality. It entails solving the Hamilton-Jacobi-Bellman-Isaacs (HJBI) Partial Differential Equation (PDE) using the method of characteristics. The hyperbolic HJBI PDE is solved using the method of characteristics with the boundary conditions exclusively specified on the Usable Part (UP) of the terminal manifold / target set which is of co-dimension 1 and where the pursuer / controller can enforce termination. The optimal state feedback control law is synthesized as opposed to obtaining an optimal control time history. The global solution is thus obtained and, in addition, the part of the state space where optimal trajectories exist is characterized. In this paper, rather than using the PMP, Isaacs’ method for the solution of Differential Games is adapted to the solution of (much simpler) optimal control problems. The importance of the UP of the terminal manifold, where the boundary conditions are specified, is emphasized. Furthermore, we derive time-optimal state feedback control laws which extend Pontryagin's canonical example concerning the regulation to a terminal state/point target to that of a target manifold of co-dimension 1. This is very much in-line with engineering practice where tolerances are specified - ``zero tolerance'' is expensive. Obviously, the size of the target sets can be shrunk to become very small. Considering non-zero tolerances, from an engineering point of view, this renders the optimal control problem easier to solve -- it obviates the need for solving a two point boundary value problem (TPBVP). Engineering and mathematics are old friends.

In this paper Isaacs' method \cite{Isaacs1965Differential} is employed rather than using the PMP to synthesize time-optimal controls for reaching a desired terminal manifold. To emphasize the advantages of using Isaacs' method rather than the PMP, we use the vehicle of the iconic example from  \cite{Pontryagin1962} which entails the dynamics of a double integrator. Time-optimal state feedback control laws are derived which globally cover the whole state space rather than constructing an optimal trajectory which leads from a specified state to a terminal state / the origin. Instead of a specific terminal state / point target, in this paper, a terminal manifold is considered. This is in tune with the engineering practice of using a finite tolerance and also renders the optimal control problem easier to solve, no TPBVP is required to be solved; if a point target must be considered the terminal manifold may be sized extremely small. We submit, that this is one way to deal with controlling to a point target in optimal control and differential games.

We consider the problem of reaching a specified target manifold, $\mathscr{C}$, of co-dimension $1$ rather than a terminal state from an arbitrary initial state in the state space in minimum time. Physical systems commonly cope with some allowable tolerance such as position or velocity error. Thus, it is the objective of this paper to investigate time-optimal control which drive a double-integrator plant which models a point mass traveling on a straight line and needs to be brought to rest at the origin in minimum time, allowing for a small error in terminal position and velocity. This also applies to the design of the roll channel autopilot of an aircraft.

Indeed, when using Isaacs' method, it becomes clear that the proper formulation of optimal control problems (and differential games) in $\mathbb{R}^n$ calls for the specification of terminal manifolds whose dimension is $n-1$. In the context of pursuit-evasion differential games, the proper treatment of ``point capture'' requires the consideration of a terminal manifold which is a sphere of radius $0<\epsilon<<1$ centered at the origin and point capture means letting $\epsilon \rightarrow 0$. This is preferable from a mathematical point of view and makes sense from and engineering point of view. And in the context of the herein discussed Pontryagin canonical example, the stipulation of a terminal manifold of co-dimension 1 brings out critical features of the optimal control problem which were hidden/obscured when ``point capture'' was considered.

The paper is organized as follows. In \Cref{sec:controlProblem}, the physical control problem is posed using non-dimensional variables. In \Cref{sec:circTargetSet}, a circular terminal manifold centered at the origin with radius $l$ is considered. The control to a non smooth target manifold with corners, a square, is investigated in \Cref{sec:squareTargetSet}. Lastly, in \Cref{sec:Conclusion} we draw conclusions.

\section{Control Problem} \label{sec:controlProblem}

Consider a point-mass with mass, $m$, which is controlled on a straight line using a bounded force, $F$. The maximum applicable force is $F_{\max}$. According to Newton's Second Law,
\begin{equation*}
    \begin{aligned}
    F(t) &= ma(t),\quad &-F_{max} \leq F(t) \leq F_{max}.\\
    \end{aligned}
\end{equation*}
Hence, the dynamics are
\begin{equation*}
    \dot{x}(t) =v(t), \quad x(0)=x_0,  \quad  \dot{v}(t) =\tfrac{1}{m}F(t), \quad  v(0) =v_0, \quad 0 \leq t\leq t_f,
\end{equation*}
where $x$ is the position on the line of the point mass, and $v$ is its velocity. The initial position of the point-mass is $x_0$ and its initial velocity is $v_0$. The goal is to drive in minimum time, $t_f$, the position and the velocity to a bounded region described as
\begin{equation}
    -L \leq x(t_f) \leq L, \quad -V \leq v(t_f) \leq V.
    \label{eq:tolSpec}
\end{equation}
It is convenient to use non-dimensional variables. The non-dimensionalization is performed as follows:
\begin{equation*}
\begin{aligned}
    x &\rightarrow x/L, & x_0 &\rightarrow x_0/L, &
    v &\rightarrow v/V, & v_0 &\rightarrow v_0/V, \\
    t &\rightarrow t\tfrac{V}{L}, & t_f &\rightarrow t_f\tfrac{V}{L}, &
    u &\triangleq \tfrac{F}{F_{max}};
    \end{aligned}
    \label{eq:nonDimForm}
\end{equation*}
where $L$ is a characteristic length and $V$ is a characteristic velocity. As is best practice in physics, also the time variable is rendered dimensionless. The dynamics in non-dimensional form are
\begin{equation*}
    \begin{aligned}
    \tfrac{dx}{dt} &= v(t), && x(0) = x_0\\
    \tfrac{dv}{dt} &= \alpha u(t), && v(0) = v_0, \quad 0 \leq t \leq t_f\\
    -1 &\leq u(t) \leq 1,
    \end{aligned}
\end{equation*}
where the non-dimensional parameter
\begin{equation*}
	\alpha \triangleq \tfrac{L F_{max}}{mV^2}.
\end{equation*}
The very same second-order dynamics/double integrator plant are encountered when designing the roll channel autopilot for an aircraft or missile; the dynamics are:
\begin{align}
	I \ddot{\varphi} = \tfrac{1}{2} \rho v^2 S b C_{l_{\delta_a}} \delta_a,
\end{align}
where $\varphi$ is the air vehicle bank angle and the control variable, $\delta_a$, is the aileron defection. $I$ is the moment of inertia, $\rho$ is the air density, $v$ is the airspeed, $S$ is the wing area and $b$ is the wing span. The non-dimensional parameter $C_{l_{\delta_a}}$ is the aircraft's lateral control derivative.
\begin{equation}
	\vert \delta_a \vert   \leq \delta_{a_\mathrm{max}}
\end{equation}
The non-dimensional dynamics of the double integrator plant are presented in \cref{eq:nonDimDyn}. Using the theory of optimal control, it is possible to design a fast lateral autopilot channel.

\begin{figure}[]
	\centering
	\begin{overpic}[width=2in]{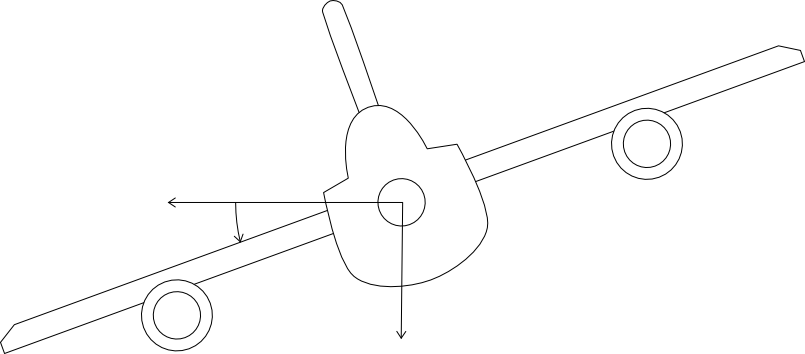}
	\put(19,21){$\hat{y}$}
	\put(51,00){$\hat{z}$}
	\put(21,14.5){$\varphi,\dot{\varphi}$}
	\end{overpic}
	\caption{Roll Autopilot Example}
\end{figure}

\section{Circular Target Set} \label{sec:circTargetSet}
\begin{figure}[]
	\centering
	\includegraphics[clip,trim=1.2in 1.2in 0in 0in,height = 2in]{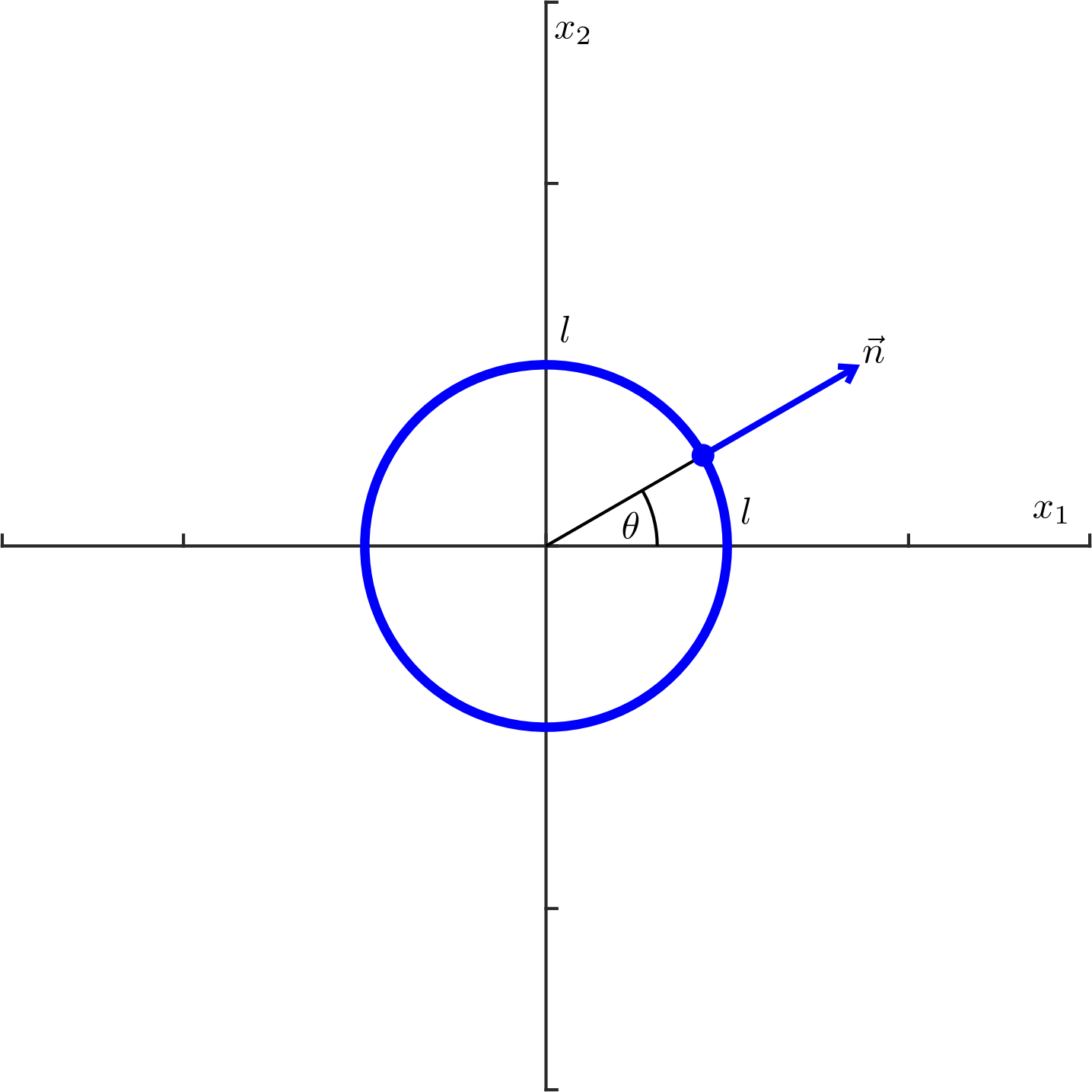}
		\caption{Circular target set with outward pointing normal $\vec{n}$}
	\label{fig:circleVectors}
\end{figure}

The physical state variables are $x_1(t) \triangleq x(t)$, $x_2(t) \triangleq v(t)$ and the non-dimensional dynamics are
\begin{equation}\begin{aligned}
     &\tfrac{d x_1}{dt} = x_2(t), x_1(0) = x_{10}, \quad \tfrac{d x_2}{dt} = \alpha u(t), x_2(0) = x_{20}, \\
     &0 \leq t \leq t_f, \quad -1 \leq u(t) \leq 1.
\end{aligned}
	\label{eq:nonDimDyn}	
\end{equation}
Consider the terminal manifold $l^2 = x^2(t_f) + \beta^2 v^2(t_f).$ The parameter $l$ is non-dimensional -- it is the tolerance parameter and $\beta$ is a non-dimensional weight parameter which trades off the importance of the terminal position error and the terminal velocity error. The optimal control problem is parameterized by $\alpha>0$, $l>0$, and $\beta>0$. For the sake of demonstration, the weight parameter $\beta$ is assumed to have the value $\beta = 1$ so we confine our attention to the terminal manifold / target set, $\mathscr{C}$, described by a circle with radius $l$ about the origin of the state space $(x_1, x_2)$ and the physical parameter $\alpha=1$.
The terminal manifold of co-dimension 1 (as required) is the circle:
\begin{equation*}
	l^2 = x_1^2(t_f) + x_2^2(t_f).
\end{equation*}
The terminal manifold of co-dimension 1 is parameterized by $0 \leq \theta \leq 2 \pi$:
\begin{equation}
	\begin{aligned}
		x_1(t_f) = l \cos \theta, \quad 
		x_2(t_f) = l \sin \theta, \quad 0 \leq \theta \leq 2\pi.
	\end{aligned}
\end{equation}
Thus, the terminal manifold
\begin{equation*}
    \mathscr{C} = \lbrace (x_1,x_2) \vert x_1 = l \cos \theta , x_2 = l \sin \theta \rbrace.
\end{equation*}
The outward pointing unit normal to the terminal manifold at $(l \cos \theta, l \sin \theta)$  is $\vec{n} = \begin{psmallmatrix}
		\cos \theta\\
		\sin \theta
	\end{psmallmatrix}$. The circular target set and an associated normal is shown in \Cref{fig:circleVectors}.

\subsection{Isaacs' Method}
The UP of the terminal manifold is where its penetration by the state can be enforced by the controller. The UP, Boundary of the Usable Part (BUP) and the Non-Usable Part (NUP) are
\begin{equation}
	\begin{aligned}
		&\text{UP} \triangleq \lbrace \mathbf{x} \vert \min_{u} \langle \vec{n}, \mathbf{f}(\mathbf{x},u) \rangle < 0 \rbrace, \quad 
		\text{BUP} \triangleq \lbrace \mathbf{x} \vert \min_{u} \langle \vec{n}, \mathbf{f}(\mathbf{x},u) \rangle = 0 \rbrace, \\
		&\text{NUP} \triangleq \lbrace \mathbf{x} \vert \min_{u} \langle \vec{n}, \mathbf{f}(\mathbf{x},u) \rangle > 0 \rbrace.
	\end{aligned}
\label{eq:Parts}
\end{equation}
For the circular terminal manifold and double-integrator system the UP, NUP and BUP -- see \cref{eq:Parts} -- are as follows
\begin{equation}
	\begin{aligned}
		\text{UP} =   \lbrace \begin{psmallmatrix} l \cos \theta \\ l \sin \theta \end{psmallmatrix}    \vert \min_{-1 \leq  u \leq 1} \left< \begin{psmallmatrix}
			\cos \theta\\
			\sin \theta
		\end{psmallmatrix}, \begin{psmallmatrix}
			x_2\\ \alpha u
		\end{psmallmatrix} \right> <0  \rbrace\\
	 \text{BUP} =   \lbrace \begin{psmallmatrix} l \cos \theta \\ l \sin \theta \end{psmallmatrix}  \vert \min_{-1 \leq  u \leq 1} \left< \begin{psmallmatrix}
	 	\cos \theta \\
	 	\sin \theta
	 \end{psmallmatrix}, \begin{psmallmatrix}
			x_2\\ \alpha u
		\end{psmallmatrix} \right> =0  \rbrace\\
	\text{NUP} =   \lbrace \begin{psmallmatrix} l \cos \theta \\ l \sin \theta \end{psmallmatrix}   \vert \min_{-1 \leq  u \leq 1} \left< \begin{psmallmatrix}
		\cos \theta\\
		\sin \theta
	\end{psmallmatrix}, \begin{psmallmatrix}
		x_2\\ \alpha u
	\end{psmallmatrix} \right> >0  \rbrace
	\end{aligned}
	\label{eq:UPBUPNUP}
\end{equation}
Having assumed, $\beta=1$, the UP will depend on the problem parameters, $\alpha$ and $l$. The UP is therefore
\begin{equation*}
	\text{UP} = \lbrace \begin{psmallmatrix} l \cos \theta \\ l \sin \theta \end{psmallmatrix} \vert \min_{-1 \leq u \leq 1} ((l \cos \theta + \alpha u) \sin \theta) < 0 \rbrace.
\end{equation*}
Therefore $\theta = 0, \theta = \pi$ are not in the UP.
\begin{enumerate}
	\item Consider the $\theta$ range $0< \theta<\pi$
	\begin{equation*}
		\text{UP}_a = \begin{cases}
			\lbrace \begin{psmallmatrix} l \cos \theta \\ l \sin \theta \end{psmallmatrix} \vert \cos^{-1} \tfrac{\alpha}{l} <\theta < \pi \rbrace & \text{if\;} \tfrac{\alpha}{l} < 1\\
			\lbrace \begin{psmallmatrix} l \cos \theta \\ l \sin \theta \end{psmallmatrix} \vert 0 < \theta < \pi \rbrace &\text{if\;} \tfrac{\alpha}{l} \geq 1		
		\end{cases}
	\end{equation*}
\item Consider the $\theta$ range $\pi< \theta<2\pi$
\begin{equation*}
	\text{UP}_b = \begin{cases}
		\lbrace \begin{psmallmatrix} l \cos \theta \\ l \sin \theta \end{psmallmatrix} \vert \pi + \cos^{-1} \tfrac{\alpha}{l} <\theta < 2 \pi \rbrace &\text{if\;} \tfrac{\alpha}{l} < 1\\
		\lbrace \begin{psmallmatrix} l \cos \theta \\ l \sin \theta \end{psmallmatrix} \vert \pi < \theta < 2\pi \rbrace &\text{if\;} \tfrac{\alpha}{l} \geq 1		
	\end{cases}
\end{equation*}
\end{enumerate}
The UP $=$ $\text{UP}_a \cup \text{UP}_b$; hence,
\begin{equation}
		\text{UP} = \begin{cases}
		\lbrace \begin{psmallmatrix} l \cos \theta \\ l \sin \theta \end{psmallmatrix} \vert 0 <\theta < \pi, \pi < \theta < 2\pi \rbrace &\text{if\;} \tfrac{l}{\alpha} \leq 1\\
		\lbrace \begin{psmallmatrix} l \cos \theta \\ l \sin \theta \end{psmallmatrix} \vert \cos^{-1} \tfrac{\alpha}{l} < \theta < \pi, \pi + \cos^{-1} \tfrac{\alpha}{l} < \theta < 2 \pi \rbrace &\text{if\;}  \tfrac{l}{\alpha} > 1		
	\end{cases}
\label{eq:UP_Circular_General}
\end{equation}

\begin{figure}[]
	\centering
	\begin{subfigure}[b]{0.48\linewidth}
		\centering
		\includegraphics[width=\textwidth]{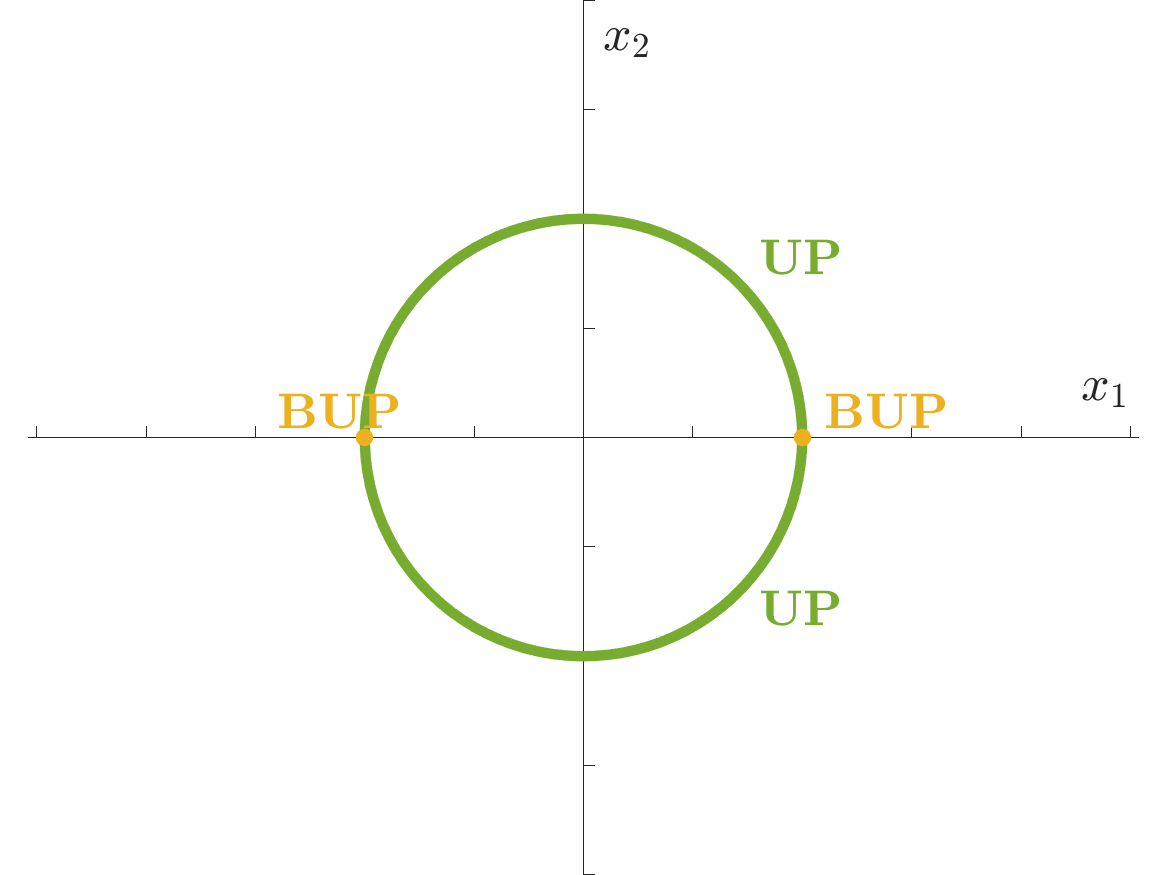}
		\caption{The circular terminal manifold\\when $0<\tfrac{l}{\alpha} \leq 1$}
		\label{fig:SmallCircle}
	\end{subfigure}
	\hfill
	\begin{subfigure}[b]{0.48\linewidth}
		\centering
		\includegraphics[width=\linewidth]{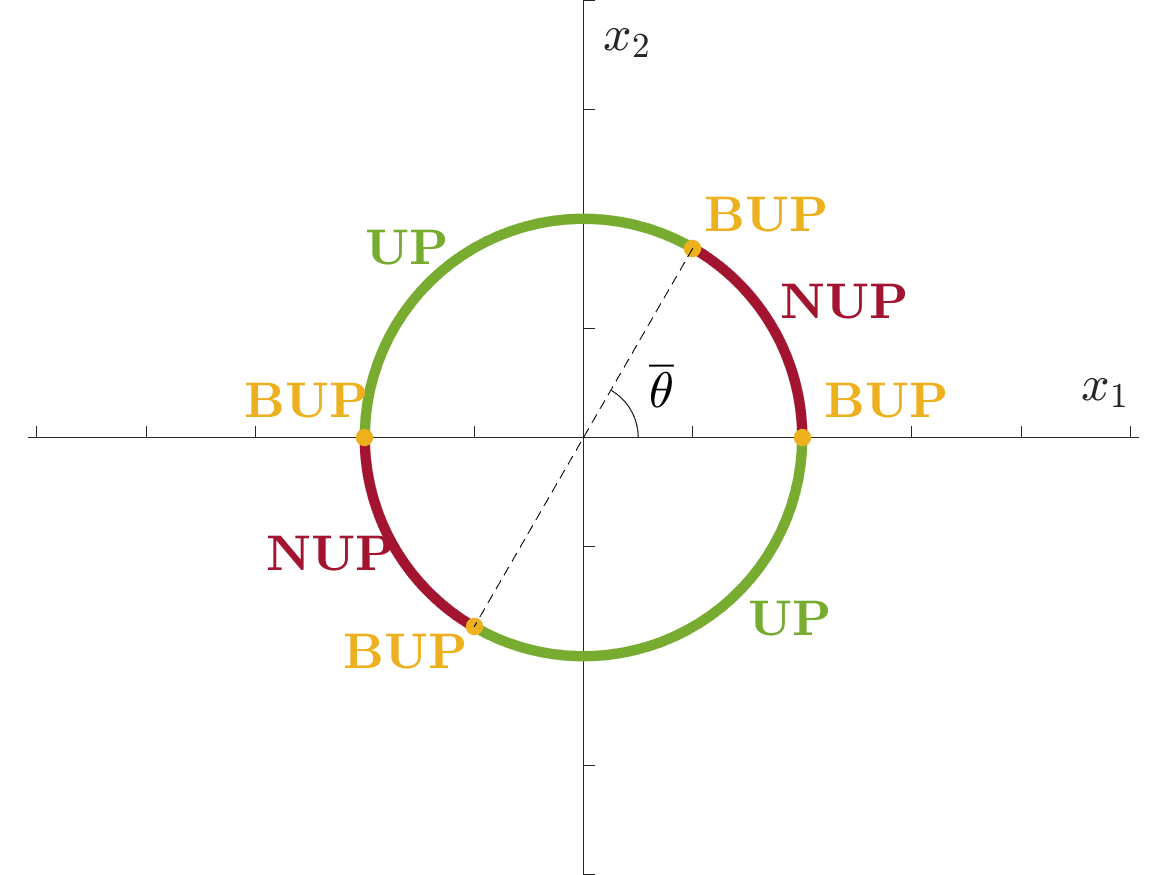}
		\caption{The circular terminal manifold\\ when $\tfrac{l}{\alpha} > 1$}
		\label{fig:LargeCircle}
	\end{subfigure}
	\caption{The UP, BUP, and NUP of the circular terminal manifold are determined by the problem parameters $l$ and $\alpha$. The NUP presents itself when $\tfrac{l}{\alpha} > 1$}
	\label{fig:CircPhaseInit}
\end{figure}

 When $\tfrac{l}{\alpha} > 1$, let $\overline{\theta} \triangleq \cos^{-1} (\tfrac{\alpha}{l})$. The UP, BUP, and NUP of the circular terminal manifold are determined by the problem parameters $l$ and $\alpha$ as described in \cref{eq:UP_Circular_General}. In \Cref{fig:CircPhaseInit}, the BUP, UP, and NUP for the circular terminal manifold is shown.
For the rest of this paper, the problem parameter is assumed to be $\alpha = 1$, so the terminal manifold depends solely upon $l$.

The Hamiltonian
\begin{equation}
	\mathscr{H} = 1 + \lambda_1(t) x_2(t) + \lambda_2(t) \alpha u(t),
	\label{eq:HamiltonianEval}
\end{equation}
with the understanding that the co-states, $\lambda_1$ and $\lambda_2$, are the partial derivatives of the Value Function with respect to the states $x_1$ and $x_2$. Dynamic programming (DP) yields the condition for optimality, $\min_u \mathscr{H}$, such that $u(t)^* = -\sign(\lambda_2(t))$. Therefore the optimal Hamiltonian is
\begin{equation*}
	\mathscr{H}^* = 1 + \lambda_1(t)  x_2(t) - \alpha \vert\lambda_2(t)\vert
\end{equation*}
$\mathscr{H}^*(t) \equiv 0, \ \forall \ 0 \leq t \leq t_f $, so
\begin{equation}
	\mathscr{H}^* \begin{large}\vert\end{large}_{t = t_f} = 0
\end{equation}

The method of characteristics employed to solve the HJBI PDE yields the Euler-Lagrage equations
\begin{equation*}
\begin{normalsize}
	\begin{aligned}
		\dot{x}_1(t) &= x_2(t), && x_1(t=0) = x_{10}\\
		\dot{x}_2(t) &= - \alpha \sign(\lambda_2(t)), && x_2(t=0) = x_{20}\\
		\dot{\lambda}_1(t) &=0, && \lambda_1(t=t_f) = a \cos \theta\\
		\dot{\lambda}_2(t) &=-\lambda_1(t), && \lambda_2(t=t_f) = a \sin \theta, && a>0
 	\end{aligned}
  \end{normalsize}
\end{equation*}
The terminal costates are also established courtesy of DP. 

We construct trajectories which emanate from the UP in \cref{eq:CircularUP}, in retrograde time $\tau>0$, and therefore we consider the retrograde dynamics,
\begin{equation}
\begin{normalsize}
	\begin{aligned}
		\mathring{x}_1(\tau) &= -x_2(\tau), && x_1(\tau=0) = l \cos \theta\\
		\mathring{x}_2(\tau) &= \alpha \text{sign}(\lambda_2(\tau)), && x_2(\tau=0) = l \sin \theta\\
		\mathring{\lambda}_1(\tau) &=0, && \lambda_1(\tau=0) = a \cos \theta\\
		\mathring{\lambda}_2(\tau) &=\lambda_1(\tau), && \lambda_2(\tau=0) = a \sin \theta,\;  \tau & \geq 0; \; \theta \in \text{UP}
 	\end{aligned}
  \end{normalsize}
 	\label{eq:RetrogradeEqns1a}
\end{equation}

The optimal Hamiltonian is zero throughout; also at the final time, and when evaluated at retrograde time, $\tau = 0$ where the co-states are specified,
\begin{equation*}
	\mathscr{H}^*\vert_{\tau = 0} = 0.
\end{equation*}
Evaluating \cref{eq:HamiltonianEval} at the final time, the coefficient $a$ is found to be:
\begin{equation}
    a = \tfrac{1}{\alpha \vert\sin \theta\vert - l \sin \theta \cos \theta},\quad \forall \; \theta \in \text{UP}
    \label{eq:aPositive}
\end{equation}
From the structure of the UP we deduce that $a$ is positive, as required. Two cases need to be considered: when $0 < \tfrac{l}{\alpha} \leq 1$ and when $\tfrac{l}{\alpha} > 1$. As will be demonstrated later the UP, BUP, and NUP as defined in \cref{eq:UPBUPNUP} differs in both cases. When, $0 < \tfrac{l}{\alpha} \leq 1$, the Usable Part (UP) and the Boundary of the Usable Part (BUP) are
\begin{equation*}
\begin{aligned}
	\text{UP} =  \lbrace \begin{psmallmatrix} l \cos \theta \\ l \sin \theta \end{psmallmatrix} \vert\;  0 < \theta < \pi,  \pi < \theta < 2\pi \rbrace, \quad 
	\text{BUP} =  \lbrace \begin{psmallmatrix} l \cos \theta \\ l \sin \theta \end{psmallmatrix} \vert\; \theta = 0, \theta = \pi   \rbrace.
 \end{aligned}
\end{equation*}
Recall, $\overline{\theta} \triangleq \cos^{-1}(\tfrac{\alpha}{l})$; when, $\tfrac{l}{\alpha} >1$,  the UP, BUP, and NUP are
\begin{equation*}
\begin{aligned}
	\text{UP} &=  \lbrace \begin{psmallmatrix} l \cos \theta \\ l \sin \theta \end{psmallmatrix} \vert\; \overline{\theta} < \theta < \pi , \overline{\theta} + \pi < \theta < 2 \pi \rbrace \\
	\text{BUP} &=  \lbrace \begin{psmallmatrix} l \cos \theta \\ l \sin \theta \end{psmallmatrix} \vert\; \theta = 0, \theta = \overline{\theta}, \theta = \pi, \theta  = \pi + \overline{\theta}  \rbrace \\
	\text{NUP} &=  \lbrace \begin{psmallmatrix} l \cos \theta \\ l \sin \theta \end{psmallmatrix} \vert\; 0 < \theta < \overline{\theta}, \pi < \theta < \pi + \overline{\theta} \rbrace.
\end{aligned}
\end{equation*}
Therefore the UP, BUP, and the NUP of the circular terminal manifold/target set are
\begin{equation}
		\label{eq:CircularUP}
		\text{UP} = \begin{cases} 
			\lbrace \begin{psmallmatrix} l \cos \theta \\ l \sin \theta \end{psmallmatrix} \vert\;  0 < \theta < \pi,  \pi < \theta < 2\pi \rbrace &\text{if\;}  0 < \tfrac{l}{\alpha} \leq 1\\
			\lbrace \begin{psmallmatrix} l \cos \theta \\ l \sin \theta \end{psmallmatrix} \vert\; \overline{\theta} < \theta < \pi , \overline{\theta} + \pi < \theta < 2 \pi \rbrace &\text{if\;} \tfrac{l}{\alpha} > 1
             \end{cases}
\end{equation}
\begin{equation}
		\label{eq:CircularBUP}
		\text{BUP} = \begin{cases} 
			\lbrace \begin{psmallmatrix} l \cos \theta \\ l \sin \theta \end{psmallmatrix} \vert\; \theta = 0, \theta = \pi   \rbrace.  &\text{if\;} 0 < \tfrac{l}{\alpha} \leq 1\\
			\lbrace \begin{psmallmatrix} l \cos \theta \\ l \sin \theta \end{psmallmatrix} \vert\; \theta = 0, \theta = \overline{\theta}, \theta = \pi, \theta  = \pi + \overline{\theta}  \rbrace &\text{if\;} \tfrac{l}{\alpha} > 1 
		\end{cases}
  \end{equation}
  \begin{equation}
		\label{eq:CircularNUP}
		\text{NUP} = \begin{cases} 
		\varnothing &\text{if\;} \tfrac{l}{\alpha} \leq 1\\
		\lbrace \begin{psmallmatrix} l \cos \theta \\ l \sin \theta \end{psmallmatrix} \vert\; 0 < \theta < \overline{\theta}, \pi < \theta < \pi + \overline{\theta} \rbrace.  &\text{if\;} \tfrac{l}{\alpha} > 1
         \end{cases}		
\end{equation}
Using the evaluation of the coefficient $a$ according to \cref{eq:aPositive} the retrograde equations in \cref{eq:RetrogradeEqns1a} are:
\begin{equation}
	\begin{aligned}
		\mathring{x}_1(\tau) &= - x_2(\tau), && x_1\vert_{\tau = 0} = l \cos \theta\\
		\mathring{x}_2(\tau) &= \alpha \text{sign}(\lambda_2(\tau)), && x_2\vert_{\tau = 0} = l \sin \theta\\
		\lambda_1 &= \tfrac{\cos \theta}{\alpha\vert\sin \theta\vert - l \sin \theta \cos \theta}\\
		\mathring{\lambda}_2(\tau) &=\lambda_1(\tau), &&  \lambda_2\vert_{\tau = 0} = \tfrac{\sin \theta}{\alpha \vert\sin \theta\vert - l \sin \theta \cos \theta}\\
		\tau &\geq 0, && \theta \in (0,\pi)\cup(\pi,2\pi)
 	\end{aligned}
 	\label{eq:RetrogradeEqns2}
\end{equation}
Therefore:
\begin{equation*}
    \lambda_2(\tau) = \tfrac{\sin \theta + \tau \cos \theta}{\alpha \vert\sin \theta\vert - l \sin \theta \cos \theta},\; \tau \geq 0,\; \theta \in (0,\pi)\cup(\pi,2\pi)
\end{equation*}

We first consider the case where the problem parameters satisfy $0<\tfrac{l}{\alpha}\leq 1$. The following abbreviations are used as necessary: $c_\theta \equiv \cos \theta$, $s_\theta \equiv \sin \theta$, $t_\theta \equiv \tan \theta$. The UP of the circular terminal manifold $\lbrace (x_1,x_2) \vert x_1^2 + x_2^2 = l^2\rbrace$ is partitioned into four quadrants as follows:
\begin{enumerate}
    \item Trajectories ``emanating'' in retrograde fashion from points on the UP of the terminal manifold which correspond to the parameter $0<\theta<\nicefrac{\pi}{2}$. For this case,  $\lambda_2(\tau)>0 \quad \forall \; \tau \geq 0$, so the optimal control $u^*(t) = -1$.
    \begin{equation*}
        \mathring{x}_2 = 1, \;   \tau \geq 0
    \end{equation*}
    Therefore
    \begin{equation*}
        \begin{aligned}
        x_2(\tau) = l \sin \theta + \tau, \quad 
        x_1(\tau) = l \cos \theta  - l \tau \sin \theta - \tfrac{1}{2} \tau ^2, \; \tau \geq 0
        \end{aligned}
    \end{equation*}
    Solving for the trajectory:
    \begin{equation*}
    \begin{aligned}
        x_1 =& - \tfrac{1}{2}x_2^2 + l \cos \theta + \tfrac{1}{2}l^2 \sin^2 \theta,\quad  x_2 > l\sin\theta, \quad \theta \in (0,\pi/2)
    \end{aligned}
    \end{equation*}

    \item Trajectories ``emanating'' from points on the terminal manifold which correspond to $\nicefrac{\pi}{2} \leq \theta<\pi$. In this case $\lambda_2(\tau)$ changes sign from positive to negative at $\tau = -\tan \theta (>0)$. For this case:
    \begin{equation*}
    \begin{normalsize}
        \mathring{x}_2 = \begin{cases} 1, &\text{if\;} 0 \leq \tau < -t_\theta \\-1, &\text{if\;} - t_\theta < \tau 
        \end{cases}
        \Rightarrow
                x_2(\tau) = \begin{cases}
        l s_\theta + \tau, &\text{if\;} 0 \leq \tau < -t_\theta\\
        l s_\theta - 2 t_\theta - \tau, &\text{if\;} - t_\theta < \tau
        \end{cases}
        \end{normalsize}
    \end{equation*}

    Therefore
    \begin{equation*}
    \begin{normalsize}
        x_1(\tau) = \begin{cases}
        l (c_\theta - \tau s_\theta) - \tfrac{1}{2}\tau ^2, &\text{if\;} 0 \leq \tau < - t_\theta\\
        \begin{aligned}
        	l (c_\theta - \tau s_\theta) + t_\theta^2 + \tfrac{1}{2}\tau^2+2 \tau t_\theta,
        \end{aligned} &\text{if\;} - t_\theta < \tau
        \end{cases}\end{normalsize}
    \end{equation*}

    Therefore when $\theta \in (\pi, \tfrac{\pi}{2})$,
    \begin{equation*}
    \begin{normalsize}
        x_1(x_2) = \begin{cases} \begin{aligned}l c_\theta- \tfrac{1}{2}x_2^2 + \tfrac{1}{2} l^2 s_\theta^2 ,\end{aligned} &\text{if\;} l s_\theta \leq x_2 < l s_\theta - t_\theta \\
        \begin{aligned}
        l c_\theta + \tfrac{1}{2}x_2^2 + 2l \tfrac{s_\theta^2}{c_\theta} - \tfrac{1}{2}l^2 s_\theta^2 - t_\theta^2 ,
        \end{aligned} &\text{if\;} l s_\theta - t_\theta  > x_2
        \end{cases}
        \end{normalsize}
    \end{equation*}

    \item Trajectories ``emanating'' from points on the terminal manifold which correspond to $\pi < \theta < \nicefrac{3\pi}{2}$. In this case $\lambda_2(\tau)$ is negative for all $\tau \geq 0$.
    \begin{equation*}
        \mathring{x}_2 = -1, \; \tau \geq 0 \Rightarrow  x_2(\tau) = l \sin \theta - \tau, \; \tau \geq 0
    \end{equation*}
    Therefore
    \begin{equation*}
        x_1(\tau) = l \cos \theta - \tau l \sin \theta + \tfrac{1}{2}\tau^2, \tau \geq 0
    \end{equation*}

    Therefore
    \begin{equation*}
        x_1(x_2) = \tfrac{1}{2}x_2^2 + l \cos \theta - \tfrac{1}{2}l^2 \sin^2\theta, \; \theta \in (\pi,\tfrac{3\pi}{2})
    \end{equation*}

    \item Trajectories ``emanating'' from points on the terminal manifold which correspond to $\nicefrac{3 \pi}{2} \leq \theta<2 \pi$. For this case $\lambda_2(\tau)$ changes sign from negative to positive at $\tau = -\tan \theta (>0)$.
    \begin{equation*}\begin{normalsize}
        \mathring{x}_2 = \begin{cases}-1,&\text{if\;} 0 \leq \tau < - t_\theta \\
        1, &\text{if\;} -t_\theta < \tau
        \end{cases}
        \end{normalsize}
    \end{equation*}
    Therefore
    \begin{equation*}
        \begin{normalsize}
        x_2(\tau) = \begin{cases}
        l s_\theta - \tau,&\text{if\;} 0 \leq \tau < - t_\theta \\
        l s_\theta + 2 t_\theta + \tau, &\text{if\;} -t_\theta < \tau
        \end{cases}
        \end{normalsize}
    \end{equation*}
    Therefore
    \begin{equation*} \begin{normalsize}
        x_1(\tau) = \begin{cases}
        l (c_\theta - \tau s_\theta) + \tfrac{1}{2}\tau^2,&\text{if\;} 0 \leq \tau < - t_ \theta \\
        \begin{aligned}&l c_\theta - t_\theta^2 -  \tfrac{1}{2}\tau ^2 - (l s_ \theta + 2 t_ \theta)\tau,
        \end{aligned}	
 &\text{if\;} -t_ \theta < \tau
        \end{cases} \end{normalsize}
    \end{equation*}
    Therefore, when $\theta \in (\tfrac{3\pi}{2},2\pi)$:
    \begin{equation*}
    \begin{normalsize}
        x_1(x_2) = \begin{cases}
        	l c_\theta + \tfrac{1}{2}x_2^2  - \tfrac{1}{2}l^2 s_ \theta^2,
         &\text{if\;} l s_\theta \geq x_2 > l s_\theta + t_ \theta \\
        \begin{aligned}
        	l c_\theta-\tfrac{1}{2}x_2^2  + 2 l \tfrac{s_\theta^2}{c_\theta} + t_\theta^2 + \tfrac{1}{2} l^2 s_ \theta^2,
        \end{aligned}
         &\text{if\;} l s_ \theta + t_ \theta > x_2
        \end{cases}
        \end{normalsize}
    \end{equation*}
\end{enumerate}

The optimal trajectories are parabolae of the form $x_1(x_2) = \pm \tfrac{1}{2} x_2^2 + c$. 
When $\theta \in (0,\tfrac{\pi}{2})\cup(\pi,\tfrac{3\pi}{2})$ no switching occurs.
When $\theta \in (\tfrac{\pi}{2}, \pi)$, $\tau_s=-\tan \theta$. We calculate $x_1(\tau_s) = \tfrac{l}{\cos\theta}-\tfrac{1}{2}\tan^2\theta$ and $x_2(\tau_s) = l\sin\theta-\tan\theta$.
When $\theta \in (\tfrac{3\pi}{2},2\pi)$, $\tau_s = -\tan \theta$. We calculate $x_1(\tau_s) = \tfrac{l}{\cos \theta} + \tfrac{1}{2}\tan^2\theta$ and $x_2(\tau_s) = l \sin \theta + \tan \theta$.

Two switching lines exist. In parametric form they are
\begin{equation}
\begin{aligned}
    x_1(\theta) = \tfrac{l}{\cos \theta} - \tfrac{1}{2}\tan^2\theta,\quad
    x_2(\theta) = l \sin \theta - \tan \theta,\; \theta \in (\tfrac{\pi}{2},\pi),
    \end{aligned}
    \label{eq:sw1}
\end{equation}
\begin{equation}
\begin{aligned}
    x_1(\theta) = \tfrac{l}{\cos \theta} + \tfrac{1}{2}\tan^2\theta,\quad
    x_2(\theta) = l \sin \theta + \tan \theta,\; \theta \in (\tfrac{3 \pi}{2},2\pi).
    \end{aligned}
    \label{eq:sw2}
\end{equation}

The switching line \cref{eq:sw1} is anchored to the circular terminal manifold at the BUP point $(-l,0)$ where $\theta = \pi$ and the switching line \cref{eq:sw2} is attached to the circular terminal manifold at the BUP point $(l,0)$ where $\theta = 0$. Define $\sigma \triangleq \tfrac{1}{\cos \theta}$. For \cref{eq:sw1}, $\sigma < 0$. Substitution of $\sigma$ in \cref{eq:sw1},
\begin{equation*}
\begin{aligned}
    &x_1(\sigma) = l \sigma - \tfrac{1}{2}\sigma^2 + \tfrac{1}{2},\quad x_2(\sigma) = (\sigma - l) \tfrac{\sqrt{\sigma^2 - 1}}{\sigma}.\\ &\sigma^2 - 2l \sigma + 2 x_1 - 1 =0
\end{aligned}
\end{equation*}
Solving for $\sigma$, subtracting $l$ from both sides, then solving for $x_2(x_1)$:
\begin{equation*}
    x_2(x_1) = \sqrt{2}\tfrac{\sqrt{l^2 + 1-2x_1}}{\sqrt{l^2 + 1-2x_1}-l}\sqrt{l^2- x_1 - l \sqrt{l^2 + 1-2x_1}}
\end{equation*}
and $x_2(-l) = 0$, as expected.\\
\begin{equation*}
    x_2(x_1) = -\sqrt{2}\tfrac{\sqrt{l^2 + 1 + 2x_1}}{\sqrt{l^2 + 1 + 2x_1}-l}\sqrt{l^2 + x_1 - l \sqrt{l^2 + 1 + 2x_1}}
\end{equation*}
and $x_2(l) = 0$, as expected.

\begin{figure}[]
	\centering
	\includegraphics[width=2.8in]{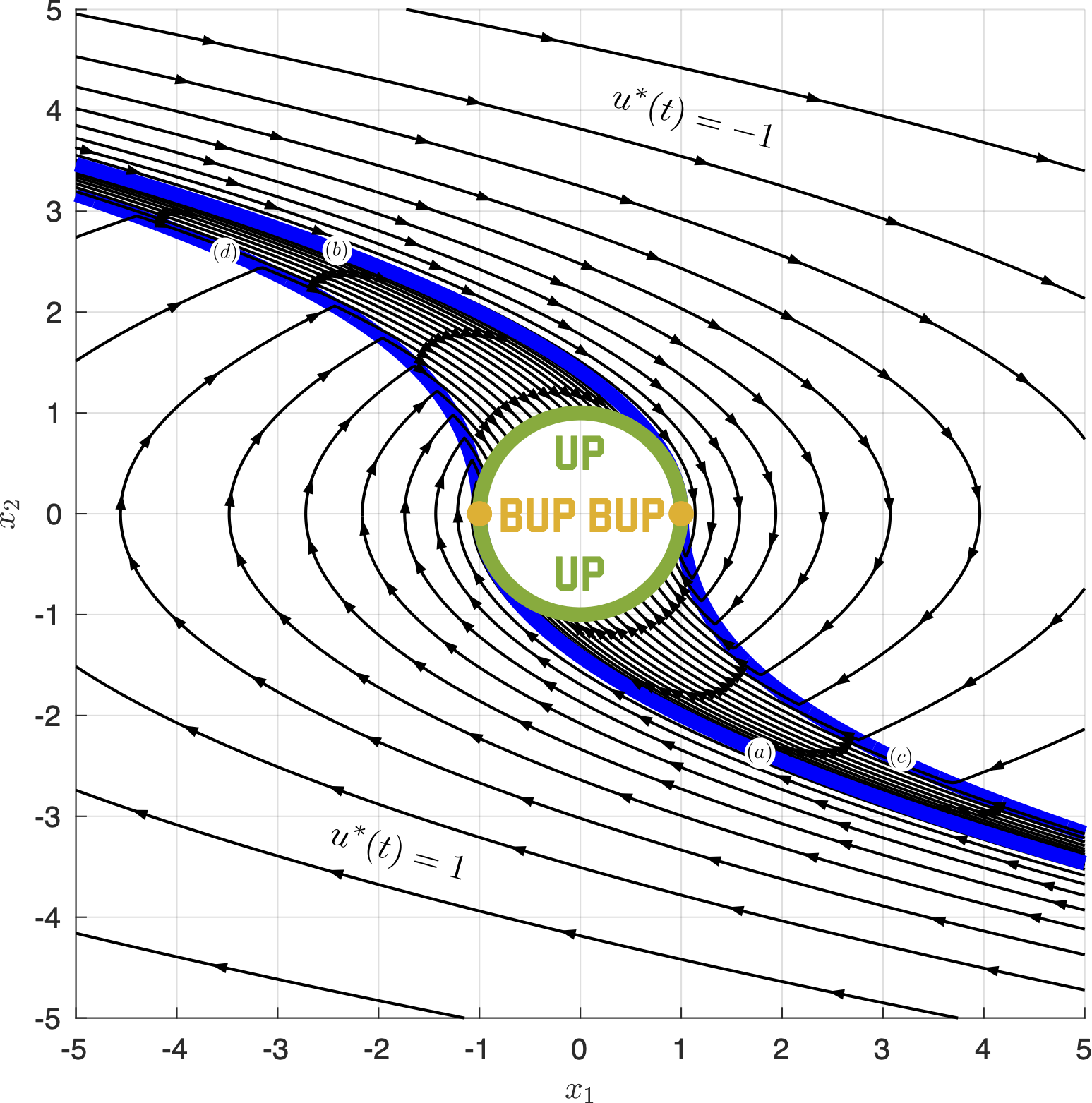}
    	\caption{The circular terminal manifold and family of trajectories in the state space when the parameters $\alpha = 1$ and $l = 1$. The curves $(a)$ and $(b)$ are where the Value Function is not continuous. At the same time, the Value Function is $C^1$. The switching lines anchored at the BUP, curves $(c)$ and $(d)$, are not optimal trajectories -- they determine when, upon crossing them, the optimal control switches from $-1$ to $1$ or from $1$ to $-1$ respectively.}
	\label{fig:circPhasel1}
\end{figure}

\begin{figure}[htbp]
	\centering
	\includegraphics[width=2.8in]{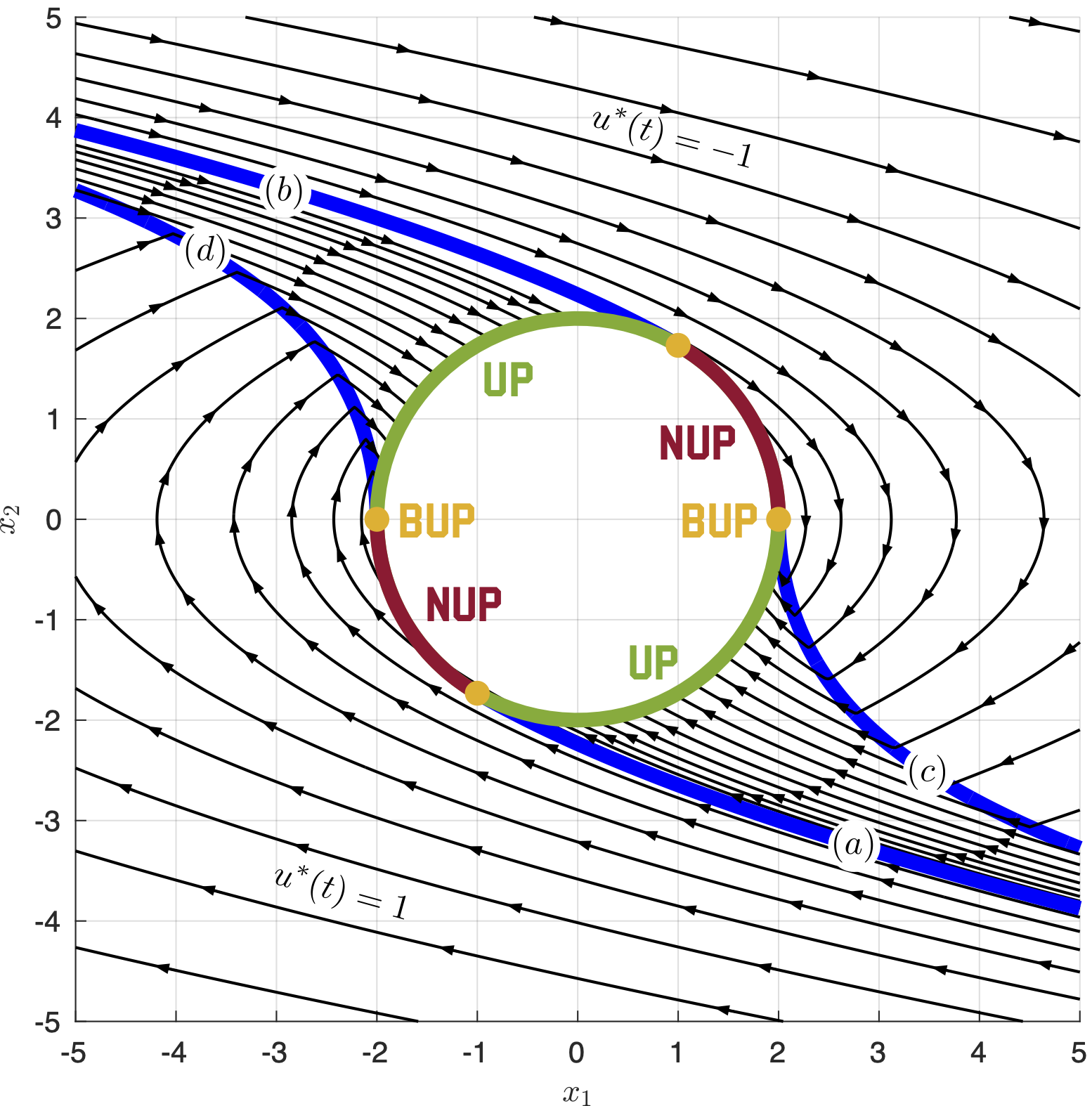}
	\caption{The circular terminal manifold and family of trajectories in the state space when the parameters $\alpha = 1$ and $l = 2$. The curves $(a)$ and $(b)$ are where the Value Function is not continuous. The switching lines, curves $(c)$ and $(d)$, determine when the optimal control switches from $-1$ to $1$ or from $1$ to $-1$ respectively. The yellow points on the boundary of the UP which do not lie on the $x_1$-axis are touch-and-go points. At these points the optimal trajectories (a) and (b) just touch the terminal manifold but do not terminate there but terminate later on the UP.}
	\label{fig:circPhasel2}
\end{figure}

The overall picture of the optimal flow field when the tolerance parameter, $0 < l \leq 1$ is shown in \Cref{fig:circPhasel1} and when the tolerance parameter $l > 1$ is shown in \Cref{fig:circPhasel2}. The parts of the terminal manifold which correspond to the UP, BUP, and NUP are indicated.
\begin{figure}[htbp]
	\begin{center}
	\includegraphics[width=3in]{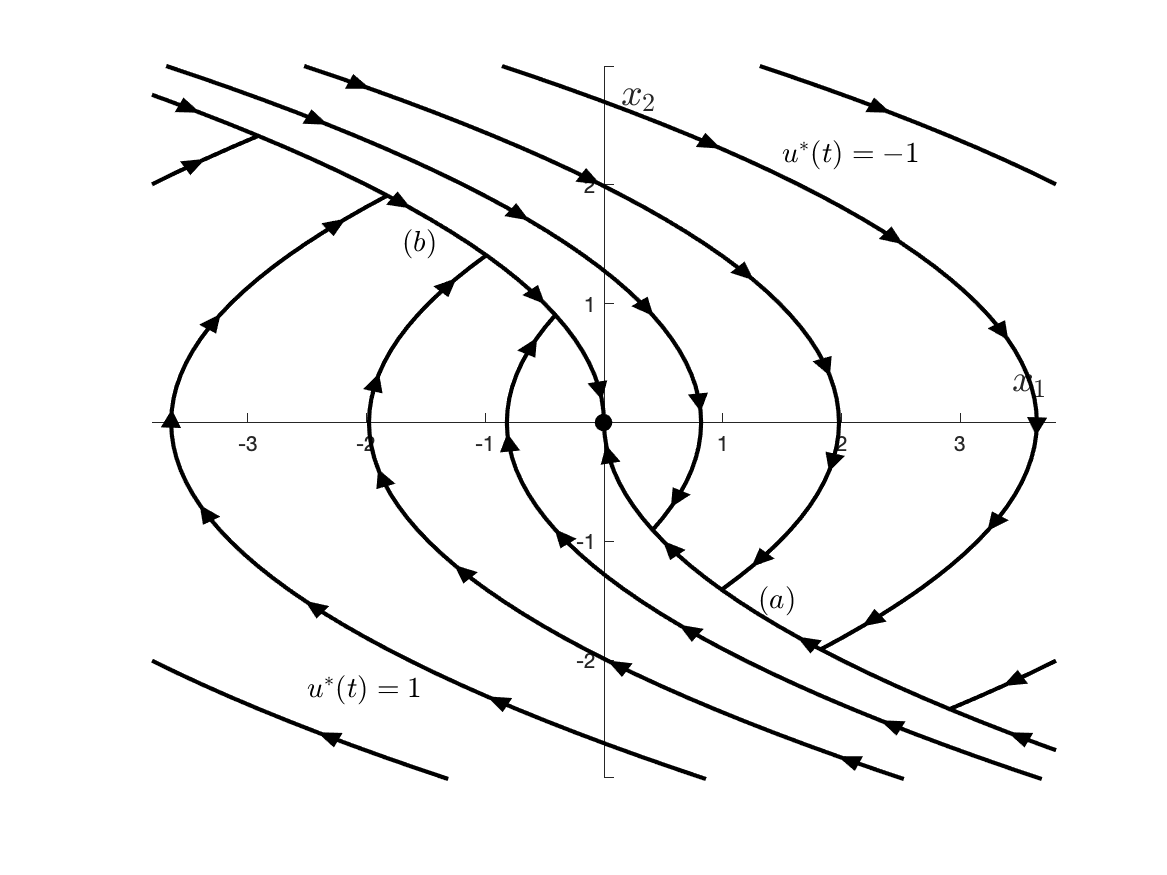}
	\caption{The ``classical'' rendition of the optimal flow field. Curves $(a)$ and $(b)$ here are optimal trajectories. Curve $(a)$  is the consolidated curve $(a)$ and $(c)$ from \Cref{fig:circPhasel1,fig:circPhasel2} and the optimal flow field trajectories contained between these two curves. Similarly, curve $(b)$ here is the consolidated curve $(b)$ and $(d)$ from \Cref{fig:circPhasel1,fig:circPhasel2} and the optimal flow field trajectories contained between these two curves. The Value function is therefore continuous.}
	\label{fig:OriginCurves}
	\end{center}
\end{figure}

In \Cref{fig:circPhasel1,fig:circPhasel2}, curves $(c)$ and $(d)$ are switching lines (SL) which are momentarily crossed by the optimal trajectories, incurring an infinitesimal loss of optimality. The SLs terminate at the BUP and therefore are not optimal trajectories themselves. Curves $(a)$ and $(b)$ delineate where the Value Function/time-to-go is not continuous. At the yellow touch-and-go points reached by curves $(a)$ and $(b)$, the trajectories do not penetrate the terminal manifold; rather, state trajectories continue on to reach the switching lines $(d)$ and $(c)$, respectively. 
The Value Function increases/jumps from above the curve $(a)$ to below curve $(a)$. The Value Function decreases/jumps from above curve $(b)$ to below curve $(b)$. Although the Value Function is not continuous, it is however, continuously differentiable away from the curves $(a)$ and $(b)$. For the sake of comparison, in \Cref{fig:OriginCurves} the canonical optimal flow field where the target set is a point target, the origin is shown. The action in \Cref{fig:circPhasel1,fig:circPhasel2} between the curves (a) and (c) and curves (b) and (d) is lost in the case of a point target because the curves $(a)$ and $(c)$ coalesce into one curve and so do curves $(b)$, $(d)$ and the optimal flow fields in-between. And these two consolidated curves have the appearance of a turnpike, that is, a universal curve in the parlance of differential games. This demonstrates the need to eschew ``point targets'' and rephrase the control problem so that it is engineering relevant. We see that \Cref{fig:OriginCurves} is quite different from the optimal flow fields shown in \Cref{fig:circPhasel1,fig:circPhasel2}. It is also interesting to present the isocost surfaces. The isocost surfaces are curves in the two-dimensional state space, $(x_1,x_2)$. The $\tau = 0$ isocost surface / curve is the UP of the target set. The case where $l=1$ is considered.
\begin{equation*}
	\text{UP} = \lbrace \begin{psmallmatrix} l \cos \theta \\ l \sin \theta \end{psmallmatrix} \vert\; 0<\theta<\pi, \pi<\theta<2\pi \rbrace
\end{equation*}

A $\tau$-isocost surface, $S_\tau$, $\tau>0$, is parameterized by $\theta$,
\begin{equation*}
	S_\tau = \Big \lbrace \begin{psmallmatrix} x_1(\theta;\tau) \\ x_2(\theta;\tau) \end{psmallmatrix} \Big \vert\; 0<\theta<\pi, \pi<\theta<2\pi \Big \rbrace
\end{equation*}

The parameter range $\lbrace \theta \vert \; 0<\theta<\pi, \pi<\theta<2\pi \rbrace$ is partitioned as follows:
\begin{equation*}
\begin{aligned}
	\lbrace \theta \; \vert& \; 0 < \theta < \pi, \pi < \theta < 2 \pi \rbrace = 
	  \lbrace \theta \;\vert\;0<\theta\leq \tfrac{\pi}{2} \rbrace  \cup \lbrace \tfrac{\pi}{2} < \theta <\pi , \tan \theta < -\tau \rbrace \\
		&\cup \lbrace \theta \vert \tfrac{\pi}{2} < \theta <\pi , \tan \theta \geq - \tau \rbrace  \cup \lbrace \theta \vert \pi < \theta \leq \tfrac{3 \pi}{2} \rbrace \\
		&\cup \lbrace \theta \vert \tfrac{3 \pi}{2} < \theta < 2\pi, \tan \theta< -\tau \rbrace \cup \lbrace \theta \vert \tfrac{3 \pi}{2} < 2\pi , \tan \theta \geq - \tau \rbrace
\end{aligned}
\end{equation*}

Let
	$\overline{\phi} \triangleq \arctan(\tau), \tau > 0$,
then,
\begin{equation*}
\begin{aligned}
\lbrace \theta\vert & 0<\theta<\pi, \pi<\theta<2\pi \rbrace =  \lbrace \theta\vert 0<\theta\leq \tfrac{\pi}{2} \rbrace \cup \lbrace \theta\vert \tfrac{\pi}{2} <\theta <\pi - \overline{\phi} \rbrace  \\
&\cup \lbrace \theta\vert \pi - \overline{\phi} \leq \theta <\pi \rbrace
    \cup \lbrace \theta\vert \pi < \theta \leq \tfrac{3 \pi}{2} \rbrace \\
 &\cup \lbrace \theta\vert \tfrac{3 \pi}{2} < \theta < 2\pi - \overline{\phi} \rbrace  \cup \lbrace \theta\vert 2\pi - \overline{\phi} \leq \theta < 2\pi \rbrace.
   \end{aligned}
\end{equation*}
The $S_\tau$ isocost surface contained within the two switching lines, $\tau>0$, is:
\begin{equation}
\begin{normalsize}
\begin{aligned}
	x_1(\theta\vert\tau) &= \begin{cases}
		l (\cos \theta - \tau \sin \theta) - \tfrac{1}{2}\tau^2 & 0 < \theta \leq \tfrac{\pi}{2}\\
		l (\cos \theta - \tau \sin \theta) - \tfrac{1}{2}\tau^2 & \tfrac{\pi}{2} < \theta < \pi - \overline{\phi}\\
		l (\cos \theta - \tau \sin \theta) + \tfrac{1}{2}\tau^2 + \tan^2\theta +2\tau \tan \theta & \pi - \overline{\phi} \leq \theta < \pi\\
		l (\cos \theta - \tau \sin \theta) + \tfrac{1}{2}\tau^2 & \pi < \theta \leq \tfrac{3\pi}{2}\\
		l (\cos \theta - \tau \sin \theta) + \tfrac{1}{2}\tau^2 & \tfrac{3\pi}{2} < \theta < 2\pi - \overline{\phi}\\
		l (\cos \theta - \tau \sin \theta) - \tfrac{1}{2}\tau^2 - \tan^2\theta - 2\tau \tan \theta & 2 \pi - \overline{\phi} \leq \theta < 2 \pi
	\end{cases}\\
	x_2(\theta\vert\tau) &= \begin{cases}
		l \sin \theta + \tau & 0 < \theta \leq \tfrac{\pi}{2}\\
		l \sin \theta + \tau & \tfrac{\pi}{2} < \theta < \pi - \overline{\phi}\\
		l \sin \theta - 2 \tan \theta - \tau \quad \quad \quad\quad \quad \quad\quad \quad \quad \; \ \ & \pi - \overline{\phi} \leq \theta < \pi\\
		l \sin \theta  - \tau & \pi < \theta \leq \tfrac{3\pi}{2}\\
		l \sin \theta - \tau & \tfrac{3\pi}{2} < \theta < 2\pi - \overline{\phi}\\
		l \sin \theta + 2 \tan \theta + \tau & 2\pi - \overline{\phi} \leq \theta < 2 \pi
	\end{cases}
\end{aligned}
\end{normalsize}
\label{eq:isocostInternal}
\end{equation}

\begin{figure}[]
	\centering
	\begin{overpic}[width=3.25in]{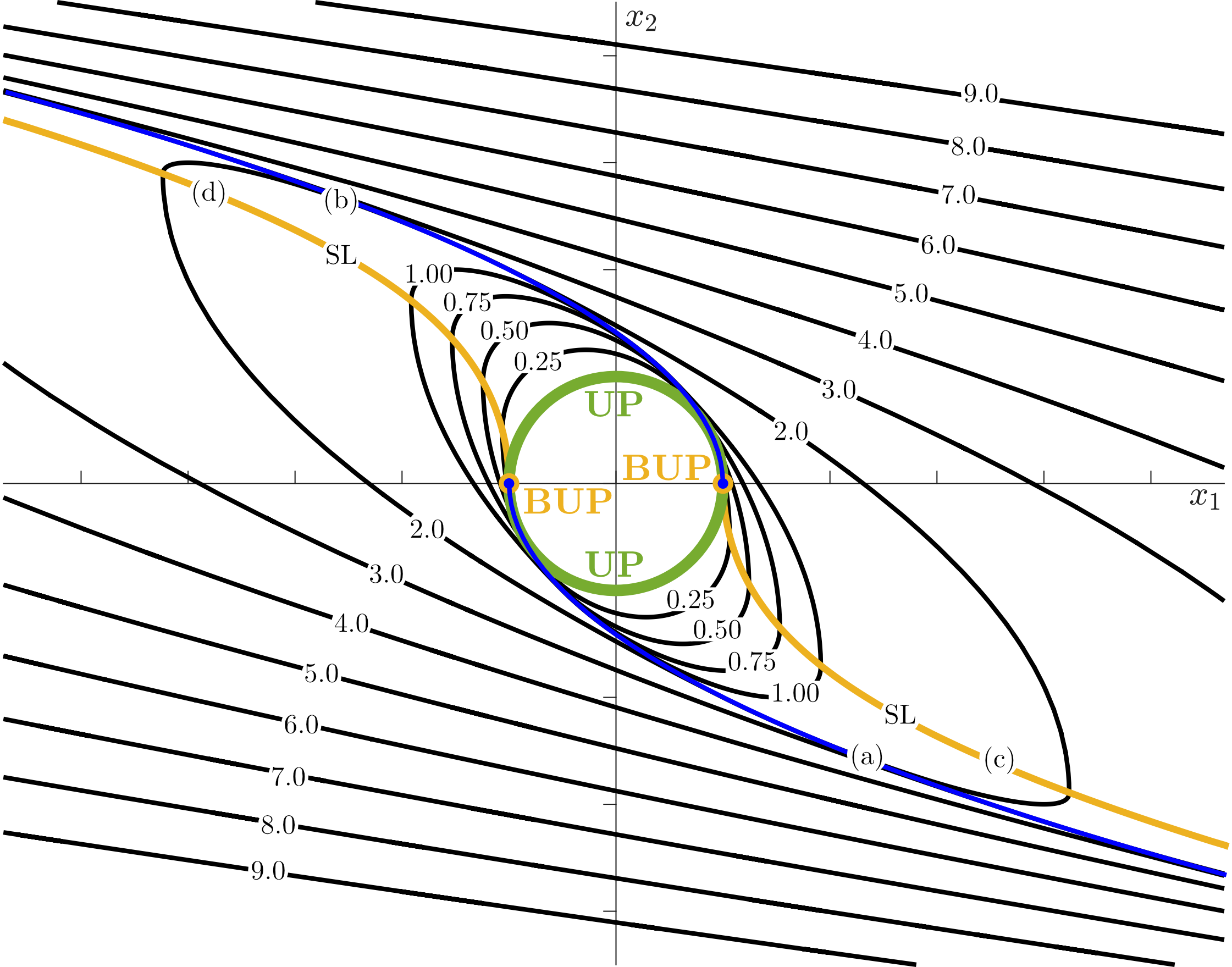}
	\end{overpic}
	\caption{The isocost/isochrone surfaces for $\tau = \lbrace 1,2, \dots, 8 \rbrace$ in the state space, $(x_1,x_2)$; $l$=1. Labeled are the curves where the value function is not continuous: (a) and (b) and the switching lines: (c) and (d).}
	\label{fig:Isocost}
\end{figure}

A figure showing the isocost curves for $\tau \in \lbrace 0.25$, $0.50$, $0.75$, $1$, $2$, $3$, $4$, $5$, $6$, $7$, $8$, $9$ $\rbrace$ are shown in \Cref{fig:Isocost}. Define the retrograde time to reach the switching line from the UP as $\tau_{1}$ and the retrograde time from the switching line to a point in the state space as $\tau_{2}$. The isocost curves of equal retrograde time $\tau$ are obtained by propogating the retrograte equations \cref{eq:RetrogradeEqns1a} back in time and constructing curves with the same propagated time, $\tau$. 

\section{Square Target Set}
\label{sec:squareTargetSet}
Recall, the tolerance specification from \cref{eq:tolSpec} repeated here for convenience: $-L \leq x(t_f) \leq L,\; -V \leq v(t_f) \leq V$.
Also let $L$ be the characteristic length. Using the non-dimensional variables as defined in \Cref{sec:controlProblem}, the terminal manifold is the square: $-1 \leq x(t_f) \leq 1,\; -1 \leq v(t_f) \leq 1.$

\begin{figure}[]
	\centering
	\begin{subfigure}[b]{0.49\linewidth}
		\centering
		\includegraphics[width=\textwidth]{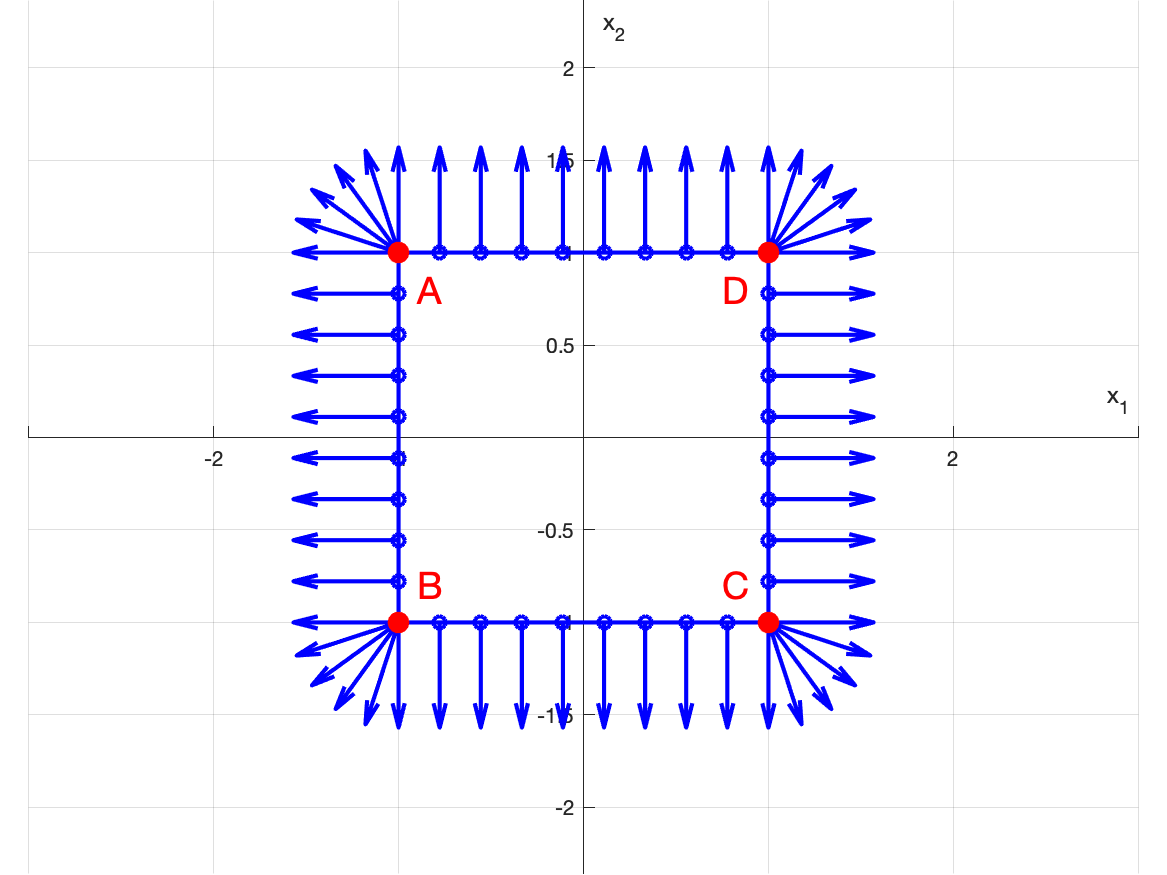}
		\caption{Outward pointing normals for the square target set $ABCD$ in the state space $(x_1,x_2)$.}
		\label{fig:squareVectorsA}
	\end{subfigure}
	\hfill
	\begin{subfigure}[b]{0.49\linewidth}
		\centering
		\includegraphics[width=\linewidth]{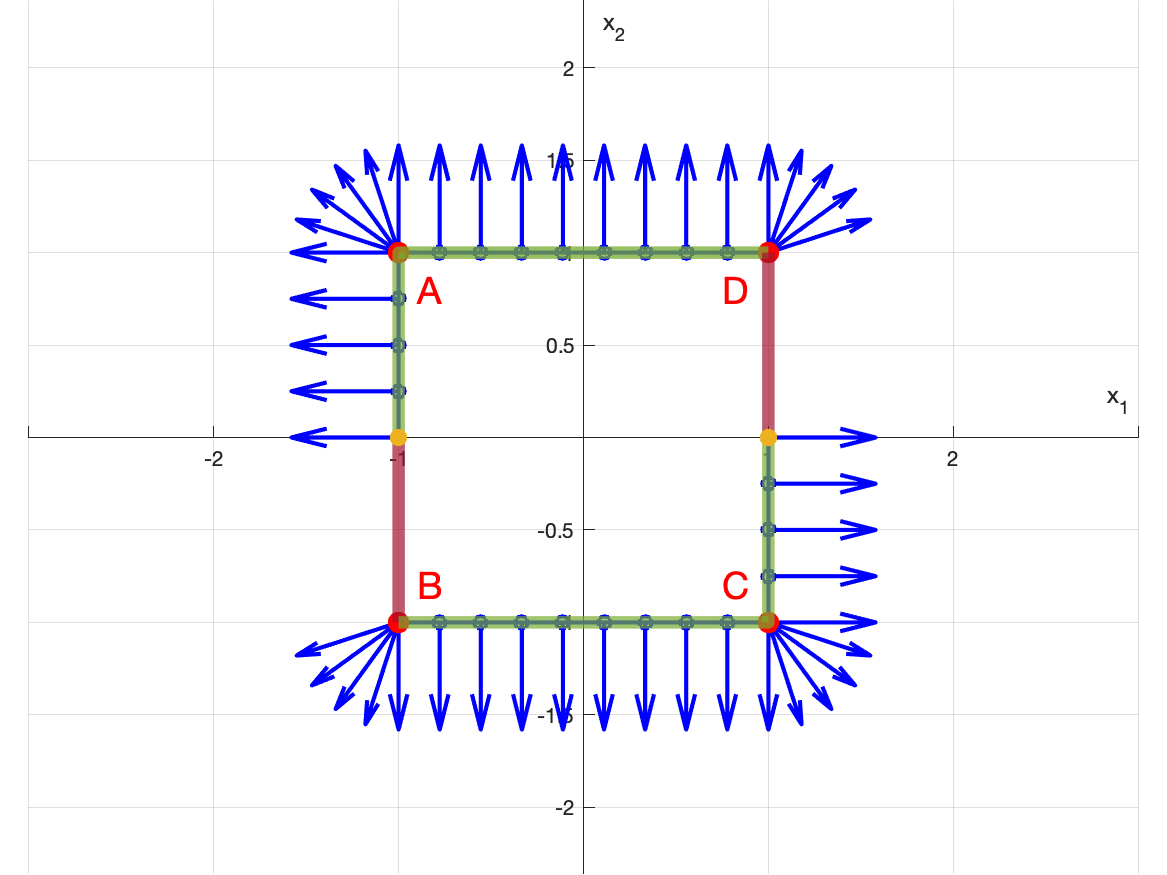}
		\caption{Outward pointing normals from the UP of the square target set $ABCD$ in the state space $(x_1,x_2)$.}
		\label{fig:squareVectorsB}
	\end{subfigure}
	\caption{The square target set $ABCD$ with outward pointing normals shown in the state space $(x_1,x_2)$. The UP of the target set is emphasized as the normals from the NUP are neglected in (b).}
	\label{fig:squareVectors}
\end{figure}

The square, non-smooth, terminal manifold $ABCD$ and the associated normals are shown in \Cref{fig:squareVectors}.
A square target set is used to show how to solve optimal control problems using Isaacs' method when the terminal manifold has corners and is not smooth. This is not apparent when considering a point target rather than a proper terminal manifold of co-dimension $1$.

\subsection{Usable Part}
In the best tradition of dynamic programs, we again ``start'' from the end. We start with identifying the UP of the terminal manifold by analyzing each of its four sides which make up the square $ABCD$, $\overline{AB}$, $\overline{BC}$, $\overline{CD}$, and $\overline{AD}$. For each side, the inner-product used to determine the UP, the BUP, and the nonusable part (NUP) are presented in \Cref{tab:UPsquareSet}. 
\begin{table}[]
    \centering
    \caption{Usable Part of the Square Terminal Manifold}
    \small
    \begin{tabular}{c   c   c  l}
        Segment & $\vec{n}$ & $\langle \vec{n} \cdot \mathbf{f} \rangle$ & Usable Part (UP)  \\   \hline  \vspace{-0.75em}
        \\ \vspace{0.75em} $\overline{AB}$ & $\begin{psmallmatrix}-1 \\0  \end{psmallmatrix}$ & $-x_2$ & $\lbrace (x_1,x_2) \vert x_1 = -1, x_2 \in (0,1] \rbrace \subset \text{UP}$ \\ \vspace{0.75em}
        $\overline{BC}$ & $\begin{psmallmatrix} 0 \\  -1\end{psmallmatrix}$ & $-u$  & $u=1 \Rightarrow \overline{BC} \subset \text{UP}$ \\ \vspace{0.75em}
        $\overline{CD}$ & $\begin{psmallmatrix} 1 \\  0\end{psmallmatrix}$  & $x_2$ & $\lbrace (x_1,x_2) \vert x_1 = 1, x_2 \in [-1,0) \rbrace \subset \text{UP}$\\ \vspace{0.75em}
        $\overline{AD}$ & $\begin{psmallmatrix} 0 \\  1\end{psmallmatrix}$ & u & $u=-1 \Rightarrow \overline{AB} \subset \text{UP}$
    \end{tabular}
    \label{tab:UPsquareSet}
\end{table}

Note that the vertices $\lbrace B,D \rbrace \notin \text{UP}$, so no optimal trajectories terminate at points $B$ and $D$. Also, at the square's vertices / corners $A$ and $C$, multiple optimal trajectories terminate because at $A$ and $C$, the normals to the terminal manifold are not unique, but form a cone whose vertex angle is $\pi/2$. And these normals are the terminal costates, each of which will give rise to an optimal trajectory. Such a family of optimal of optimal trajectories will contribute to forming the optimal flow field which must cover the entire state space. The terminal state is

\begin{equation}
    \begin{normalsize}
    x_1(t_f) = \begin{cases}
    -1   &\text{on} $\;\;$ \overline{AB} \\
    s_1  , -1 < s_1 \leq 1  &\text{on} $\;\;$ \overline{BC} \\
    1     &\text{on} $\;\;$ \overline{CD} \\
    s_2  , -1 \leq s_2 < 1  &\text{on} $\;\;$ \overline{AD}
    \end{cases},
    \end{normalsize}
    \ 
    \begin{normalsize}
    x_2(t_f) = \begin{cases}
    s_3 , 0 < s_1 \leq 1  &\text{on} $\;\;$ \overline{AB} \\
    -1    &\text{on} $\;\;$ \overline{BC} \\
    s_4 , -1 \leq s_4 < 0    &\text{on} $\;\;$ \overline{CD} \\
    1   &\text{on} $\;\;$ \overline{AD}
    \end{cases}
    \end{normalsize}
    \label{eq:xtf}
\end{equation}
The terminal co-states are aligned with the outward pointing normals
\begin{equation}
    \begin{normalsize}
\begin{aligned}
    &\lambda_1(t_f) = \begin{cases}
    -a_1 & \text{on} $\;\;$ \overline{AB} \\
    0    & \text{on} $\;\;$ \overline{BC} \\
    a_2  & \text{on} $\;\;$ \overline{CD} \\
    0    & \text{on} $\;\;$ \overline{AD} \\
    a_5 \cos \theta_1 & \text{at} $\;\;$ A, \theta_1 \in (\tfrac{\pi}{2},\pi) \\
    a_5 \cos \theta_2 & \text{at} $\;\;$ C, \theta_2 \in (\tfrac{3 \pi}{2}, 2\pi)
    \end{cases}\\
    &a_1>0,\; a_2 > 0,\; a_5 > 0
\\
    &\lambda_2(t_f) = \begin{cases}
    0    & \text{on} $\;\;$ \overline{AB} \\
    -a_3 & \text{on} $\;\;$ \overline{BC} \\
    0,   & \text{on} $\;\;$ \overline{CD} \\
    a_4  & \text{on} $\;\;$ \overline{AD} \\
    a_6 \sin \theta_1 & \text{at} $\;\;$ A, \theta_1 \in (\tfrac{\pi}{2},\pi) \\
    a_6 \sin \theta_2 & \text{at} $\;\;$ C, \theta_2 \in (\tfrac{3 \pi}{2}, 2\pi)
    \end{cases}\\
    &a_3>0,\; a_4 > 0,\; a_6 > 0
\end{aligned}
\end{normalsize}
\label{eq:lambdatf}
\end{equation}

Therefore a family of optimal trajectories will terminate at points $A$ and $C$. 

Recall the Hamiltonian, $\mathscr{H} = 1 + \lambda_1 x_2 + \lambda_2 u$, and since $u^* = -\text{sign}(\lambda_2)$, the optimal Hamiltonian, $\mathscr{H}^* = 1 + \lambda_1 x_2 - \text{sign}(\lambda_2)\lambda_2$. The optimal Hamiltonian is zero, including at the final time, just as described in \Cref{sec:circTargetSet}. The coefficient $a$ is determined by evaluating the Hamiltonian at the final time, $t_f$ where the co-states are known. This is accomplished by substitution of the values from \cref{eq:xtf,,eq:lambdatf} into \cref{eq:HamiltonianEval}. The resulting values for $a_1$ through $a_6$ are
\begin{equation}
\begin{aligned}
    a_1 &= \tfrac{1}{s_3}, 0 < s_3 \leq 1 &a_2 &= -\tfrac{1}{s_4},  -1 \leq s_4 < 0, \quad a_3 = 1, \quad a_4 = 1,\\
    a_5 &= \tfrac{1}{\sin \theta_1 - \cos \theta_1},  \tfrac{\pi}{2} \leq \theta_1 \leq \pi,
    &a_6 &= \tfrac{1}{\cos \theta_2 - \sin \theta_2},  \tfrac{3 \pi}{2} \leq \theta_2 \leq 2 \pi
\end{aligned}
    \label{eq:aParameters}
\end{equation}
Note that $a_5>0$ and $a_6>0$ over the domain of $\theta_1$ and $\theta_2$ as required. Moreover, substitution of the $a$ parameters from \cref{eq:aParameters} into the co-state equations \cref{eq:lambdatf} provides the final co-states.
\begin{equation}
\begin{aligned}
    &\begin{normalsize}
    \lambda_1(t_f) = \begin{cases}
    -\tfrac{1}{s_3},\; 0 < s_3 \leq 1 & \text{on} \;\; \overline{AB} \\
    0    & \text{on} \;\; \overline{BC} \\
    -\tfrac{1}{s_4},\; -1 \leq s_4 < 0  & \text{on} \;\; \overline{CD} \\
    0    & \text{on} \;\; \overline{AD} \\
    \tfrac{\cos \theta_1}{\sin \theta_1 - \cos \theta_1} & \text{at} \;\; A \\
    \tfrac{\cos \theta_2}{\cos \theta_2 - \sin \theta_2} & \text{at} \;\; C
    \end{cases}
\end{normalsize}, \
\begin{normalsize}
	\lambda_2(t_f) = \begin{cases}
    0    & \text{on} \;\; \overline{AB} \\
    -1 & \text{on} \;\; \overline{BC} \\
    0,   & \text{on} \;\; \overline{CD} \\
    1  & \text{on} \;\; \overline{AD} \\
    \tfrac{\sin \theta_1}{\sin \theta_1 - \cos \theta_1} & \text{at} \;\; A \\
    \tfrac{\sin \theta_2}{\cos \theta_2 - \sin \theta_2} & \text{at} \;\; C
    \end{cases}
\end{normalsize},\\
&\text{where} \ \tfrac{\pi}{2} \leq \theta_1 \leq \pi, \tfrac{3 \pi}{2} \leq \theta_2 \leq 2\pi 
\end{aligned}
\label{eq:lambdatfEval}
\end{equation}

The Euler-Lagrage / characteristic equations are
\begin{equation*} 
	\begin{aligned}
	\dot{x}_1(t) &= x_2(t),  &x_1(t=0) &= x_{10}, \\ \dot{x}_2(t) &= -\sign(\lambda_2(t)),  &x_2(t=0) &= x_{20},\\
	\dot{\lambda}_1(t) &=0,  &\lambda_1(t=t_f) &= \text{Eq. \cref{eq:lambdatfEval}}, 
		\\ \dot{\lambda}_2(t) &=-\lambda_1(t),  &\lambda_2(t=t_f) &= \text{Eq.  \cref{eq:lambdatfEval}}
 	\end{aligned}
\end{equation*}
In retrograde time, $\tau$, we consider trajectories which emanate from the UP, and therefore we have
\begin{equation}
	\begin{aligned}
		\mathring{x}_1(\tau) &= x_2(\tau), &x_1(\tau=0) &= \text{Eq. \cref{eq:xtf}}, 
		\\\mathring{x}_2(\tau) &= -\text{sign}(\lambda_2(\tau)), &x_2(\tau=0) &= \text{Eq. \cref{eq:xtf}}, \\ 
		\mathring{\lambda}_1(\tau) &=0, &\lambda_1(\tau=0) &= \text{Eq. \cref{eq:lambdatfEval}},
		\\\mathring{\lambda}_2(\tau) &=-\lambda_1(\tau), &\lambda_2(\tau=0) &= \text{Eq.  \cref{eq:lambdatfEval}}, \; 
		 \tau \geq 0
 	\end{aligned}
 	\label{eq:RetrogradeEqns}
\end{equation}
Because the derivative for $\lambda_1$ is zero, $\lambda_1(\tau) = \lambda_1(\tau =0), \tau \geq 0$. Using this information, we calculate, in retrograde, $\lambda_2(\tau),\; \tau \geq 0$. Integrating $\mathring{\lambda}_2(\tau)$
\begin{equation*}
\begin{normalsize}
	\begin{aligned}
		&\lambda_2(\tau) = \begin{cases}
			-\tfrac{1}{s_3}\tau,\; 0<s_3\leq1 & \text{on} \;\; \overline{AB} \\
			-1& \text{on} \;\; \overline{BC}  \\
			-\tfrac{1}{s_4}\tau,\; -1\leq s_4<0 & \text{on} \;\; \overline{CD} \\
			1 & \text{on} \;\; \overline{AD} \\
			\tfrac{\tau \cos \theta_1 + \sin \theta_1}{\sin \theta_1 - \cos \theta_1 },\; \tfrac{\pi}{2}\leq \theta_1 \leq \pi & \text{at} \;\; A\\
			\tfrac{\tau \cos \theta_2 + \sin \theta_2}{\cos \theta_2 - \sin \theta_2},\; \tfrac{3 \pi}{2} \leq \theta_2 \leq 2 \pi & \text{at}\;\; C
		\end{cases} \qquad \tau \geq 0
	\end{aligned}
 \end{normalsize}
\end{equation*}
Next, we calculate the optimal trajectories in retrograde fashion.
\begin{equation}
    \begin{normalsize}
        \begin{aligned}
            &(x_1(\tau),x_2(\tau)) = \begin{cases}
		(-1 -s_3\tau +\tfrac{\tau^2}{2},\;s_3 -\tau)& \text{on} \;\; \overline{AB}\\
		(s_1 + \tau +\tfrac{\tau^2}{2},\;-1 -\tau) & \text{on} \;\; \overline{BC}\\
		(1-s_4\tau - \tfrac{\tau^2}{2},\; s_4+\tau) & \text{on} \;\; \overline{CD}\\
		(s_2 - \tau-\tfrac{\tau^2}{2}, \; 1+\tau) & \text{on} \;\; \overline{AD}\\
		(-1 -\tau -\tfrac{\tau^2}{2},\;1 + \tau), \;0\leq\tau\leq  -\tan \theta_1 & \text{at} \;\; A\\
		\begin{aligned}(-1 - \tau + 2 \tau \tan \theta_1  + \tfrac{\tau^2}{2} + \tan^2\theta_1, \\ \; 1 - 2 \tan \theta_1 - \tau) ,\quad  -\tan\theta_1 \leq \tau \end{aligned} & \text{at} \;\; A\\
		(1 +\tau + \tfrac{\tau^2}{2},\;-1 -\tau), \;0\leq\tau\leq  -\tan \theta_2 & \text{at} \;\; C\\
		\begin{aligned}(1 + \tau - 2 \tau \tan \theta_2 - \tfrac{\tau^2}{2}-\tan^2\theta_2
		,\\-1 +2 \tan \theta_2 + \tau), -\tan\theta_2 \leq \tau\end{aligned} &\text{at} \;\; C
	\end{cases}\\
        &\qquad \quad -1 < s_1 \leq 1, \quad -1 \leq s_2 <1 \\
	&\qquad \quad0 < s_3 \leq 1, \quad -1\leq s_4 <0 \\
	&\qquad \quad \tfrac{\pi}{2} \leq \theta_1 \leq \pi, \tfrac{3 \pi}{2} \leq \theta_2 \leq 2\pi
        \end{aligned}
    \end{normalsize}
\end{equation}
The optimal trajectories with $x_1$ as a function of $x_2$, are
\begin{equation}
\begin{normalsize}
\begin{aligned}
	&x_1(x_2) = \begin{cases}
		\tfrac{x_2^2-s_3^2-2}{2}, & x_2 \leq s_3, \mathbf{x}_f \in \overline{AB}\\
		\tfrac{x_2^2-1+2s_1}{2}, & x_2\leq -1,  \mathbf{x}_f \in \overline{BC}\\
		\tfrac{s_4^2-x_2^2+2}{2}, & x_2 \geq s_4,   \mathbf{x}_f \in \overline{CD}\\ 
		\tfrac{1-x_2^2+2s_2}{2}, &  x_2\geq 1,  \mathbf{x}_f \in \overline{AD}\\ 
		\tfrac{-x_2^2-1}{2}, & \begin{aligned}&1\leq x_2 \leq 1 - \tan \theta_1,            
            \mathbf{x}_f \in A\end{aligned}\\
		\tfrac{x_2^2-3}{2} + 2\tan \theta_1 - \tan^2\theta_1,
		&  \begin{aligned}
  	&x_2 \leq 1 - \tan \theta_1, \mathbf{x}_f \in	A
  \end{aligned}\\
		\tfrac{x_2^2+1}{2}, & \begin{aligned} &\tan\theta_2 - 1 \leq x_2 \leq -1, \mathbf{x}_f \in C	
 \end{aligned}\\
		\tfrac{-x_2^2+3}{2} - 2\tan \theta_2 + \tan^2\theta_2,
		  & \begin{aligned}
 	&\tan \theta_2 -1 \leq x_2, \mathbf{x}_f \in	C
 \end{aligned}\\
	\end{cases}\\
	&\qquad \quad -1 < s_1 \leq 1, \quad -1 \leq s_2 <1 \\
	&\qquad \quad0 < s_3 \leq 1, \quad -1\leq s_4 <0 \\
	&\qquad \quad\tfrac{\pi}{2} \leq \theta_1 \leq \pi, \tfrac{3 \pi}{2} \leq \theta_2 \leq 2\pi
\end{aligned}
\end{normalsize}
\label{eq:squarex2x1}
\end{equation}

\begin{figure}[]
	\centering
	\begin{overpic}[percent,width=3in]{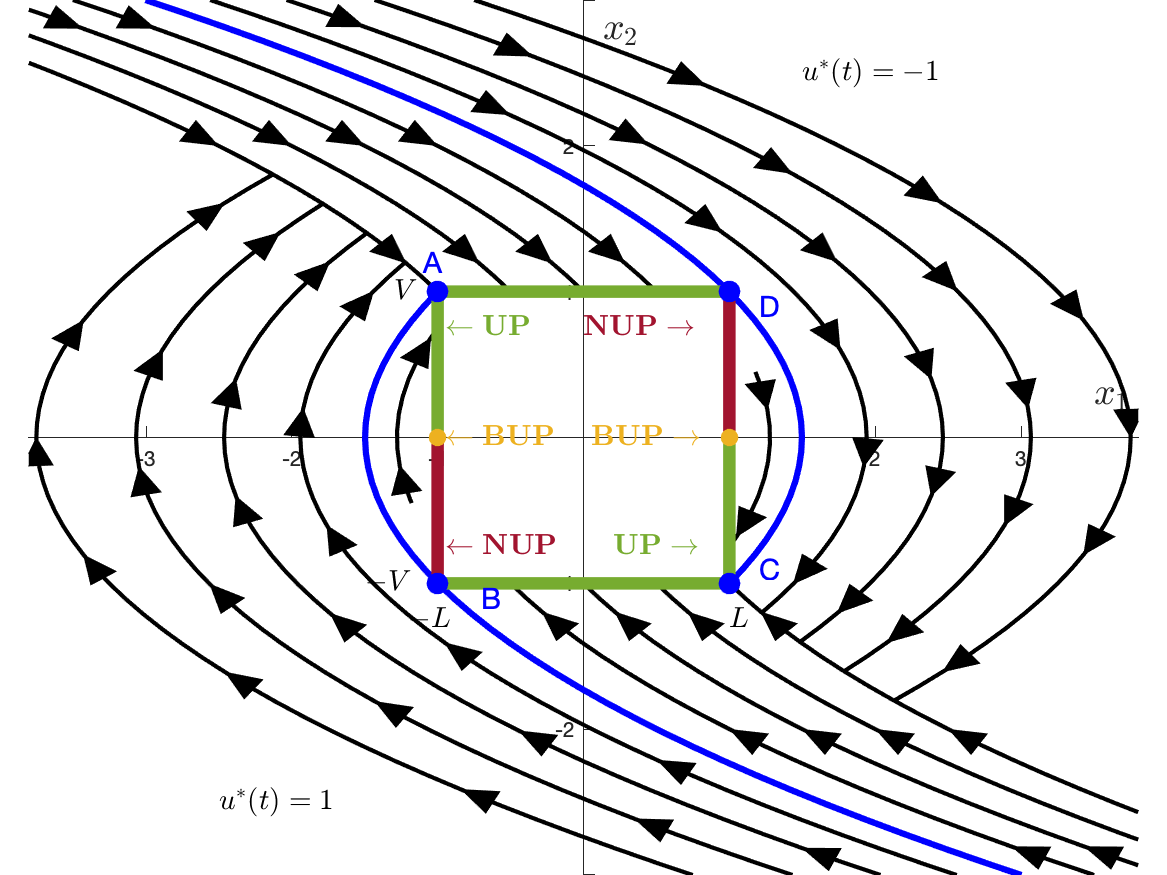}
	\put(8,62){$(a)$}
	\put(85,12){$(b)$}
	\end{overpic}
	\caption{The optimal flow-field and two switching curves, (a) and (b), for the square target manifold.}
	\label{fig:squarePhasePortrait}
\end{figure}

Notice that the sign of $\lambda_2$ changes at the time instant $\tau_s = -\tan \theta_1$ for trajectories emanating from point A, and $\tau_s = -\tan \theta_2$ for trajectories emanating from point C. The trajectories specified in \cref{eq:squarex2x1} provide the optimal trajectories for $x_1$ and $x_2$, provided a single parameter: $s_1$, $s_2$, $s_3$, $s_4$, $\theta_1$, or $\theta_2$. These trajectories potentially fill the two-dimensional state space: $(x_1,x_2)$.
\noindent Two switching curves exist:
\begin{enumerate}
    \item The switching curve pertaining to the family of optimal trajectories which terminate at vertex, $A$ of the terminal manifold (and which are parameterized by $\pi/2 \leq \theta_1 \leq \pi$ is:
	\begin{equation}
		x_1 = -\tfrac{1}{2} x_2^2 - \tfrac{1}{2},\quad x_2 \geq 1
		\label{eq:swPointA}
	\end{equation}
    \item The switching curve pertaining to the family of optimal trajectories which terminate at the vertex $C$ of the terminal manifold (and which are parameterized by $3\pi/2 \leq \theta_2 \leq 2\pi$) is:
	\begin{equation}
		x_1 = \tfrac{1}{2}x_2^2 + \tfrac{1}{2}, \quad x_2 \leq -1
		\label{eq:swPointC}
	\end{equation}
\end{enumerate}
No switching occurs on the four families of optimal trajectories which terminate on the four sides of the terminal manifold. The optimal flow field for reaching the square target manifold is shown in \Cref{fig:squarePhasePortrait}.

In \Cref{fig:squarePhasePortrait} the blue lines represent switching lines. The switching line which is anchored at point $A$ is described by \cref{eq:swPointA}; while the switching line which is anchored at point $C$ is described by \cref{eq:swPointC}. The switching lines (a) and (b) are, themselves, optimal trajectories for reaching the UP since they are anchored at point $A$ and $C$ respectively. The switching lines (a) and (b) are when the optimal control used for reaching the UP of the terminal manifold switches from 1 to -1 and from -1 to 1, respectively.

The corners $B$ and $D$ are not in the UP and therefore trajectories which pass thought these points continue and terminate at points $A$ and $C$ respectively. These are ``touch-and-go'' trajectories. The touch-and-go trajectory which passes though point $B$ has the equation: $x_1 = x_2^2/2 - 3/2$ and the touch-and-go trajectory which passes through point $D$ has equation $x_1 = -x_2^2/2 + 3/2$. The square target set might be small in extent; yet, all of this is missed when the target set is a point, as in Pontryagin's canonical example.

For the sake of visualization, consider the mental exercise whereby a small circular arc of radius $\epsilon <<1$ is located at the corner of our square  terminal manifold, rendering it a smooth terminal manifold. A unique normal is then associated with each point on this small circular arc, and consequently, a unique optimal trajectory terminates at this point on the  small circular arc. Consequently, a family of optimal trajectories terminates on this small circular arc, and by extension, at the corner of the manifold. We have literally smoothed the rough edges of a terminal manifold. And in the parlance of optimal control theory, we have put to  work the viscosity solution concept.

\section{Conclusion} \label{sec:Conclusion}
In this paper, the PMP and DP methods for solution of optimal control problems are juxtaposed. We advocate solving optimal control problems by leveraging Isaacs' constructive method for the solution of differential games. Isaacs' method is based on the method of DP as opposed to the PMP which is rooted in the calculus of variations and provides necessary conditions for optimality, which however can afford the construction of an optimal trajectory. We also emphasize the importance of formulating well-posed optimal control problems in an engineering way, that is, the need to move away from ``point capture'' and instead consider terminal manifolds of co-dimension $1$; point capture is then the limiting case where the terminal manifold is shrunk to a singleton and the terminal manifold is a point target. This is aligned with engineering practice where finite tolerances are specified.

Most importantly, in the paper we draw attention to the fact that in the process of correctly solving an optimal control problem, identifying the  UP of the terminal manifold is \emph{the} critical first step. It comes down to properly specifying the boundary conditions of the HJBI PDE for the Value Function. Thus, in minimum time optimal control problems the Value Function is correctly set to zero, we emphasize, \emph{only} on the UP of the terminal manifold. Blindly setting the Value Function to  zero on the whole terminal manifold just because the state is in the terminal manifold, and proceeding with the numerical solution of the HJBI PDE, yields trajectories which are not either feasible or optimal. All of this is observed when point targets are considered; that's why one can say that when a point target is considered, the problem is not well posed from an engineering point of view.

Two examples are used to highlight the solution of min-time optimal control problems for reaching both a smooth and non-smooth terminal manifold, rather than a point target, where this distinction is obscured. This not only renders the optimal control problem easier to solve, but from an engineering point of view, is realistic, in that it represents an acceptable terminal tolerance/error. The first case is a circular terminal manifold, the second is a square terminal manifold with corners. The former highlights how to pose and solve optimal control problems when the terminal manifold is smooth, while the latter highlights the solution process when the terminal manifold has corners. While Isaacs' method naturally requires the terminal manifold being of co-dimension 1, it becomes apparent the classical PMP based approach hides critical aspects of the optimal control problem when a point target is considered; this is highlighted in this paper. While the PMP is a necessary condition of optimality, which however allows the construction of candidate optimal trajectories, Isaacs' method directly yields the global optimal flow field state feedback control laws, and as a byproduct, also the region of controllability.

\section*{Data Availability}
Data sharing not applicable to this article as no datasets were generated or analyzed during the current study.

\end{document}